\documentclass{article}
\usepackage{amsmath}
\usepackage{subcaption}
\usepackage{multicol}
\usepackage{multirow}
\usepackage{PRIMEarxiv}

\usepackage[utf8]{inputenc} 
\usepackage[T1]{fontenc}    
\usepackage{hyperref}       
\usepackage{url}            
\usepackage{booktabs}       
\usepackage{amsfonts}       
\usepackage{nicefrac}       
\usepackage{microtype}      
\usepackage{lipsum}
\usepackage{fancyhdr}       
\usepackage{graphicx}       
\graphicspath{{media/}}     

\title{A bi-objective optimization model to plan vaccination campaigns
}

\author{
  Fernando Montenegro-Dos Santos \\
  School of Engineering \\
  Universidad Católica del Norte \\
  Coquimbo, Chile\\
  \texttt{fernando.montenegro@ce.ucn.cl} \\
   \And
  Francisco Pérez-Galarce \\
  Department of Computer Science \\
  Pontificia Universidad Católica de Chile \\
  Santiago, Chile\\
  \texttt{frperezga@gmail.com} \\
 \AND
   Carlos Monardes-Concha \\
   School of Engineering \\
  Universidad Católica del Norte \\
  Coquimbo, Chile\\
   \texttt{cmonardes@ucn.cl} \\
   \And
   Sergio Cruz-Zárate \\
   Pukará Research and Consulting \\
   Santiago, Chile \\
   \texttt{scruz@pukara-research.cl} \\
   \And
   Alfredo Candia-Véjar \\
   School of Engineering \\
   Universidad Finis Terrae \\
    Santiago, Chile\\
   \texttt{alfredo candia.vejar@gmail.com} \\
}

\begin{document}
\maketitle

\begin{abstract}
Vaccination campaigns have saved thousands of lives, reaching the farthest places in the world. These campaigns have required substantial investments and accurate coordination between several actors within the vaccine supply chain. Despite these successful strategies, the outbreak of COVID-19 has altered the objectives and rules of undertaking vaccine campaigns. Then, it is essential to consider two new facts in planning vaccination campaigns. First, some groups of infected people by the virus are more vulnerable to severe illness. Second, the virus is exceptionally contagious; sometimes, no symptoms are apparent. Accordingly, we propose a bi-objective optimization model that allows healthcare decision-makers to design effective vaccination campaigns by considering these COVID-19 characteristics and controlling the associated costs. Careful utilization of temporary and traditional vaccination centers is crucial to creating a more robust strategy. Located in suitable places, temporary centers help increase the probability of reaching groups difficult to be vaccinated while simultaneously avoiding crowd congestion and reducing the risk of spreading infections in dispensing vaccination centers. Experiments were conducted using data from an area in Santiago, Chile. The results show the model allows us to manage the resources by providing a variety of vaccination plans according to the prioritization level of each objective.

\end{abstract}


\section{Introduction and Motivation}
A large number of pandemics have occurred in human history. In ancient times, one of them was the Plague of Athens (430 to 426 BC) in which typhoid fever killed a quarter of the population during the Peloponnesian War. The Black Death was a bubonic plague pandemic which occurred in Afro-Eurasia during the fourteenth century. It was the most fatal pandemic in recorded history, killing between 75 and 200 million people in Europe, Asia, and North Africa. The 1918-1920 Spanish flu infected half a billion people around the world, including on remote Pacific islands and in the Arctic, killing 20 to 100 million. Influenza outbreaks in 1918 had a high mortality rate for young adults. 
Coronaviruses (CoV) are a large family of viruses that cause illnesses ranging from the common cold to more severe diseases such as Middle East Respiratory Syndrome (MERS-CoV) and Severe Acute Respiratory Syndrome (SARS-CoV-1). 

A new strain of coronavirus (SARS-CoV-2) caused Coronavirus disease 2019, or COVID-19, which was declared a pandemic by the World Health Organization (WHO) on 11 March 2020. As of October $4^{th}$, 2022, more than 6.5 million people have died \cite{WorldHealthOrganization2022WHODashboard}. COVID-19 has expanded throughout the world, and after one year and a half,  several strategies, such as social distancing and vaccination, have been partially effective to contain the pandemic \cite{cadeddu2022planning}. Vaccination campaigns are being applied worldwide, and preliminary data shows that can help to protect against the virus. The WHO introduced a practical phasing for managing epidemics and pandemics. The phases are Anticipation, Early detection, Containment, Control and Mitigation, and Elimination or Eradication \cite{WorldHealthOrganization2018ManagingDisease}. As a matter of fact, several logistics problems will have to be solved. In the phase of Control and Mitigation, one key problem is the planning of vaccination campaigns. In this paper, we focused on the  "last-mile delivery" that refers to the dispensing of the vaccine to the people.

To reach susceptible groups, the effective design of large-scale vaccination campaigns worldwide is a major challenge for many governments and non-governmental entities \cite{cadeddu2022planning, perrone2021influenza}. To achieve a significant vaccination impact on the target population, it is essential to have an efficient vaccine supply chain and accurate planning \cite{Lee2017TheWorld}. For planning large-scale vaccination campaigns, several aspects must be taken into account: (i) we must focus on specific priority target groups, such as health teams, vulnerable populations, adults, children, or people from six months to 18 years; (ii) we must identify the places where the different groups will be vaccinated; (iii) we need to consider the epidemiological behavior of the disease to define the vaccination period; (iv) we should act as soon as possible and in a short time, according to the availability of resources; and finally, (v) technical concerns associated with the coordination between agencies, budgets and communication must be addressed \cite{PanAmericanHealthOrganization2009RegionalPandemic}.

Typically, vaccination campaigns can be considered standard health operations in which long-term and mid-term planning are suitable. However, when extreme events emerge, some outbreaks pose complex logistic scenarios similar to those challenges observed during humanitarian logistic contexts (e.g. earthquakes, tsunamis or hurricanes). Humanitarian Logistics (HL) considers a well-know cycle composed of the following four stages: \textit{mitigation}, \textit{preparedness}, \textit{response} and \textit{recovery}; these phases have been discussed and described in \cite{Altay2006OR/MSManagement, Tomasini2009HumanitarianLogistics, Galindo2013ReviewManagement}.
There is a fundamental difference between the life cycle in HL and that of a pandemic. Traditional disasters are determined by one or two disruptive events (for example an earthquake followed by a tsunami) that have a short duration (hours or a few days) and normally affect a large region. So the phases in HL are defined as being either pre-disaster or post-disaster operations. 

A pandemic is a process that is disseminated throughout a large region or even around the world, like COVID-19. This process is complex because it typically takes two or more years and sometimes remains active for decades (a current example is HIV). The WHO introduced a practical phasing for managing epidemics and pandemics. The phases are Anticipation, Early detection, Containment, Control and mitigation, and Elimination or eradication \cite{WorldHealthOrganization2018ManagingDisease}. Our model is associated with the \textit{Control and mitigation} phase of the pandemic. In this paper, we provide a bi-objective optimization model to provide a vaccine distribution plan. This model focuses on the last mile of a large and complex supply chain described below. In this context, the "last mile delivery" refers to the dispensing of the vaccine to the people. Therefore, the different phases of HL are intermixed in a pandemic, rendering the problem more complex. 

A classic vaccine supply chain (VSC) considers different factors and stages, such as manufacturers, distribution centers (DCs), transport, points of dispensing, and temporary centers. Throughout this network, different actors must make many decisions under uncertainty to complete the vaccination process \cite{Duijzer2018LiteratureChain}. Additionally, within VSCs, it is also possible to differentiate between the part of the network that moves vaccines from one country to another and the part that moves vaccines within a country. As previously commented, this paper focuses on the distribution network in the last mile phase of dispensing the vaccine to the population, specifically in an urban context. 

In this context, despite Operations Research (OR) tools have been used widely to address different topics in the design of VSCs (e.g., the selection of facilities to serve as vaccination dispensing points; inventory issues such as the periodic review and updating of supplies; transportation and  the distribution of supplies; and the purchase of supplies), there is scarce literature about the last mile problem that occur once vaccines arrive at the health centers where the vaccines are normally applied \cite{Farahani2019ORDevelopments}.

Beyond ad-hoc mathematical models, several lessons and new rules have been learned from a series of outbreaks in the last decades (e.g., H1N1-2009 and yellow fever in Brazil in 2018). Furthermore, the ongoing vaccination campaign for SARS-CoV-2 exposes the logistic complexity in a massive distribution with a latent contagious risk in the environment and the outbreak evolution \cite{Teytelman2012ModelingPopulation}. From experience gained in the last decades, we can highlight the following beneficial practices: (i) the use of new channels (e.g., drive-through, schools, nursing homes, and temporary centers) to distribute the vaccinations has been a good source of innovation among practitioners to improve vaccination campaigns' efficiency and coverage; (ii) the segmentation of target groups by the risk level has been widely used to minimize the risk in more vulnerable groups. Obviously, each practical lesson forces us to rethink our optimization models to increase their usability. Therefore, we propose an optimization model that incorporates some of these lessons.

In the context of COVID-19, new critical logistics problems appeared during 2020 in main health centers. The installed infrastructure in health centers for specialized beds to infected people by the virus was surpassed. Thus, other rooms assigned to other diseases had to be reassigned to treat infected people by COVID-19. Also, the allocation of ventilators was another problem faced by the decision-makers in the health system.  Research about logistics problems related to the management of vaccination campaigns have been recently published. \cite{Grundel2021HowOutbreaks} focused on the planning horizon on the optimal vaccination/social distancing strategy depending on the long run or too short-horizon. \cite{Molla2021AdaptiveGroups} used optimal control theory and heuristic strategies for allocating vaccines against COVID-19. They formulated a mathematical model that tracks disease spread among different age groups and across different geographical regions. \cite{Bertsimas2021WhereFacilities} studied the problem of vaccination prioritization to mitigate the near-end impact of the pandemic. It is important to comment that these papers did not study tactical/operative plans for vaccination campaigns. For example, the definition of places for vaccine dispensing in an environment of high infection rates is not considered.

In short-term and mid-term, it is expected that one or more vaccine campaigns will be necessary and will be applied to a significant portion of the population of (almost) all ages. 

Modeling risks in the vaccination process can be challenging since, as we have observed in the ongoing global campaign, numerous factors are involved, such as supplier reliability, risk of each target group, duration of the campaign, stock of vaccines, availability of vaccination personnel, and availability of permanent and temporary facilities for vaccination. But mainly, in the context of a virus transmitted by close contact, it is mandatory to consider congestion in the vaccination centers.

The following contributions are presented in this paper with the aim to decrease the gap research concerning logistics problems in vaccination campaigns for COVID-19.
\begin{itemize}
    \item We propose a new model-based policy for planning the last-mile single-dose vaccine dispensing, which considers the complex environment concerning COVID-19. A combined modeling of permanent and temporary vaccination places that provides more protection and prioritizes certain groups of people are designed to follow some basic characteristics from humanitarian logistics. The problem is formulated as a bi-objective optimization model that minimizes the risk associated with congestion, balanced with a second objective of minimizing the costs associated with temporary centers. The set covering model locates temporary vaccination centers and assigns groups of people to both types of centers (temporary and permanent), respecting the given priorities of the campaign.

    \item To evaluate the model's behavior in a case study, we used data from   the San Bernardo commune in Santiago, Chile. The model results showed that effective scheduling of the different groups was done in the horizon planning by the location of a group of temporary centers. Additionally, a prioritization of the groups in the time was also achieved. From these experiences, some managerial insights are discussed in the case study.
    
   \item  We provide a set of managerial insights to simplify the application of our model. For example, we highlight the impact of some crucial parameters on the vaccination campaign, such as the maximum number of temporary center or the weights of each objective function.
\end{itemize}

This article continuous as follows. Section \ref{section2} presents a literature review of mathematical models for planning vaccination campaigns. Then, in Section \ref{section3}, the operational problem is described, indicating actors, objectives, and decisions. Then, a bi-objective formulation for this problem and an illustrative example is presented. Next, in Section \ref{section4}, the case study is presented with the data collection and related assumptions. In Section \ref{section5}, experiments from several viewpoints are conducted, and a set of managerial insights are given. Finally, conclusions and future work are also presented in Section \ref{section6}.

\section{Literature Review}
\label{section2}

The literature review is divided three-fold. Section \ref{section2.1} concerns discussing recent papers focused on logistical problems associated with COVID-19 vaccination. Then, Section \ref{section2.2} is about review papers that analyze approaches to solve logistical problems in the context of other virus caused-diseases. Finally, Section \ref{section2.3} explores some specific logistic issues in vaccination campaigns connected with our situation.

\subsection{Logistics problems with COVID-19 vaccination}
\label{section2.1}
After the worldwide expansion of the COVID-19 and its terrible consequences on human beings since 2020, one main challenge has been the development of vaccines to help to control the effects of the virus. Thus, scientific literature focused on reporting advances in the impact of the various vaccines. In 2021, vaccine campaigns are being applied all over the world. Lately, some research reports have been studying relations between mobility, social distancing, and the effect of vaccines. These studies have implications on the logistics operations to be carried out to the vaccination campaigns, like prioritizing the groups to be vaccinated. Therefore, we first focused on the discussion about these papers.

\cite{Grundel2021HowOutbreaks} studied the impact of the planning horizon on the optimal vaccination/social distancing strategy from the design of a compartmental model that extends a previous model of this type. They found that leading to reduce social distancing in the long run, and it is essential to vaccinate the people with the highest contact rates first. They commented that this step should not overburden the healthcare system. It is important to note from this study the relevance of having a flexible vaccination planning system in the dynamic process of the COVID-19. They also found that if the objective is to minimize deaths, the same priority should be applied, provided that the social distancing measures are sufficiently strict. Furthermore, they also showed that if too short-horizon are used, only the elderly are vaccinated, and more social distancing is necessary.

\cite{Molla2021AdaptiveGroups} used optimal control theory and heuristic strategies for allocating vaccines against COVID-19. They formulated a mathematical model that tracks disease spread among different age groups and across different geographical regions. They also proposed a method to combine age-specific contact data with geographical movement data. Their method was applied in mainland Finland utilizing mobility data from a major telecom operator. The tested scenarios found that distributing vaccines demographically and in an age-descending order is not optimal for minimizing deaths and disease burden. Instead, more lives could potentially be saved by strategies emphasizing high-incidence regions and distributing vaccines parallel to multiple age groups. They noted the importance of updating the vaccination strategy when the effective reproduction number changes due to the general contact patterns changing and new virus variants entering.  

\cite{Bertsimas2021WhereFacilities} studied the problem of vaccination prioritization to mitigate the near-end impact of the pandemic. They integrate an epidemiological model to capture the effects of vaccinations and the variability in mortality rates across age groups into a prescriptive model to optimize the location of vaccination sites and subsequent vaccine allocation.
The problem is formulated as a bilinear, non-convex optimization model. They proposed a coordinate descent algorithm that iterates between optimizing vaccine allocations and simulating the dynamics of the pandemic. The model analyzes allocating vaccines in the "spots" of the pandemic versus administering vaccines to the most vulnerable risk groups. They commented that experimental results using real-world data in the United States suggest that the proposed optimization approach can yield more benefits than a benchmark allocation that distributes vaccines proportionally to each subpopulation size. No reference to the operative problem of where the population would be vaccinated is done in the paper.

\cite{Georgiadis2021OptimalChain} proposed an optimization model for the optimal planning of the COVID-19 vaccine supply chain. The model tries to find optimal decisions regarding vaccine operations between locations, the inventory profiles of central hubs and vaccination centers, and the daily vaccination plans in the vaccination centers of the supply chain network. The model was tested on a study case that simulates the Greek nationwide vaccination program.

\subsection{Logistical problems in vaccination operations}
\label{section2.2}
This section is focused on discussing the contribution of some review papers concerning the logistics problems in vaccination. \cite{Jordan2021OptimizationReview} reviewed the literature on optimization models and methods in the context of COVID-19. The review focus is on formal optimization techniques or machine learning approaches. In particular, in the section on Mass Vaccination, they found an article \cite{Matrajt2021VaccineFirst} where an age-stratified model was used in an optimal vaccine allocation plan based on three metrics: deaths, symptomatic infections, and hospitalizations. They also found an article \cite{Jadidi2021AData} considering the allocation prioritization of certain groups (e.g., age, presence of comorbidities) and strategies for vaccine distribution. 
 
\cite{Dasaklis2012EpidemicsReview} discussed the importance of logistics operations for the control of epidemic outbreaks. In their analysis of the logistics network design, the authors signaled main aspects such as defining the location, quantity, and capacity of the storage centers and points of dispensing, allocating these facilities to serve specific people groups, and establishing transport and inventory policies. The authors included various mathematical programming models and optimization algorithms in the problem's modeling and solution. \cite{Dasaklis2012EpidemicsReview}

\cite{Duijzer2018LiteratureChain} emphasized the benefits of developing a supply chain perspective instead of partially analyzing logistics problems. They claimed that in the VSC were coordination problems as many parties with conflicting interests participated in the process. Additionally, they found a difference between developed and developing countries in the distribution stage mainly due to a lack of quality infrastructure.

\cite{DeBoeck2020VaccineReview} reviewed vaccine distribution chains in low and middle-income countries by two sides: (1) characteristics and challenges of the VSC, and (2) seven classification criteria for the reviewed OR papers, including the modeled part of the VSC and the countries' and vaccines' coverage. The analysis concluded that most papers focus more on the strategic decision level than the tactical or operational. Though several papers indicate various performance measures, they do not show how to prioritize these measures. There is an opportunity to use multi-criteria decision-making in VSC management. Keeping this in mind, our proposal model, situated on the tactical and operative levels, incorporates a multi-objective function that balances the risk of contagion per priority group and the cost of using temporary vaccination centers.

\cite{tang2022bi} proposed a bi-objective MIP model to plan a vaccination campaign. The model seeks to minimize, in the first place, the total operational cost of the campaign and, second, minimize the total population distance to vaccination sites. The decisions considered are: 1) How many sites to open, 2) How many servers per site to open, 3) How many doses to replenish per site and day, and 4) How to assign people to each site to get vaccinated. To solve this model, they used two algorithms: the weighted-sum method and the $\epsilon$-constraint method. To solve a practical example instance, they also proposed a genetic algorithm. The model lack considering the contagion risk between risk groups of people. Therefore no congestion aspects were considered in the model.

\subsection{Logistics operations in vaccination campaigns}
\label{section2.3}
Finally, the following papers discussed logistics and operations related to the vaccine campaigns before the COVID-19. Still, some of these operations are common to the current problem in the context of COVID-19, like last-mile distribution, uncertainties in the vaccine supply network, and facility location.

\cite{Lemmens2016AChains} claimed that vaccines are different from commodity goods, and the VSC needs special considerations. Specifically, for network design, it is relevant to consider the following characteristics: (1) allocation, (2) location, (3) limited shelf life and cold chain distribution, (4) production capacity planning, and (5) batch sizing. \cite{Lin2020ColdChain} discussed the logistic problem related to cold chain distribution. \cite{tsai2022solving} focused on patient allocation to face epidemic Dengue fever. \cite{Enayati2020OptimalEquity} proposed a mathematical programming model to optimize the distribution of influenza vaccines. The model minimizes the number of vaccine doses distributed to faced the early stages of the pandemic. They demonstrated that using subgroups based on geographic location and age on the model is the key factor in controlling the transmission dynamic for optimal distribution. \par

\cite{Lim2022RedesignNetworks} proposed a mixed-integer programming (MIP) model to design an Expanded Program on Immunization (EPI) VSC and a hybrid approach that combines metaheuristic and exact methods to solve it. The model decides the location of a set of DCs, the flow paths from the country's central store through one or more health clinics where vaccination occurs, the amount of the dose in each DC, the vehicles, and the number of trips required of each vaccine delivery. \cite{delgado2022equity} focused on the importance of deciding on the opening of DCs, based on the fact that the location of these facilities must guarantee equitable distribution and the centers must be accessible and meet the corresponding demand. They propose a nonlinear programming model to solve this problem based on a coverage formula that uses the Gini index to measure equity and accessibility.

\cite{Sadjadi2019TheApproach} presented a multi-objective optimization model for designing the vaccine supply network under uncertain conditions. They focused on the network's strategic and tactical aspects, considering the following objectives: minimizing total cost and unsatisfied high- and low-priority demands. They used the weighting sum method to find its solution. Also, \cite{Azadi2018StochasticProducts} developed a stochastic optimization model to define vaccine vials' distribution strategies in developing countries. The decisions were modeled to be made weekly, and the uncertainty of the weekly amount of patient arrivals was considered to improve inventory replenishment cycles.
\cite{Rachaniotis2012AStudy} studied the logistic problem related to optimizing the use of scarce resources to face a national vaccine plan for influenza. The model seeks to minimize the total number of infections by influenza, determining the prioritization of subpopulations to be vaccinated. The model uses mobile medical teams to reach targeted people who cannot move to local vaccination centers.

\cite{Yang2021OptimizingCountries} proposed an improved vaccine distribution network, which is an alternative to the traditional structure, in which vaccines could be distributed using a more flexible network that can select facilities appropriately from the legacy network. They used a MIP and proposed a novel algorithm that solves a sequence of increasingly large-size MIP problems. \cite{Yang2021OutreachApproach} presented an outreach approach to improving the distribution of vaccines in low and middle-income countries. The idea behind this approach was to set up mobile clinics at distant locations to move the vaccination process closer to people living in rural or remote areas. They modeled the problem as a MIP that combined a set covering problem with a vehicle routing problem. They also considered uncertainty concerning the population (and hence the volume of vaccines required) at each location and the travel times.

In conclusion, several papers have discussed multiple aspects of the VSC, formulating optimization models and applying mathematical programming and (meta)heuristics methods. To our best knowledge, no research has addressed vaccine campaign planning in a complex environment like the current pandemic. COVID-19 has generated some new protocols for coexistence in society, such as social distancing. This component is critical in health systems, so a new strategy to implement the dispensing of vaccines to the population is proposed in the next section.

\section{An optimization model for planning vaccination campaigns} 

\label{section3}

This section presents a bi-objective optimization model for planning single-dose vaccination campaigns, considering temporary and permanent facilities. Assumptions and related definitions are outlined in Section \ref{assumptions}. In Section \ref{mathematical_formulation}, the mathematical formulation for this problem is shown. Finally, in Section \ref{3.3}, some solutions and insights from this model are given using an illustrative example. 

\subsection{Basic concepts and assumptions}
\label{assumptions}

Our problem assumes that the vaccine distribution chain's last mile occurs when vaccines arrive at hospitals, clinics, or health centers. After that, the population is vaccinated in the previously mentioned facilities.  However, the new scenario of vaccination during the ongoing COVID-19 pandemic needs the utilization of a combination of different facilities. One of these facilities is the usual vaccination center (permanent centers). In some cities, several of these centers exist. However, due to the necessity of social distancing, they cannot vaccinate many people. Therefore, it is necessary to define new appropriate facilities that can supplement the capacity of the permanent centers. We call them temporary centers in a generic designation, such as gyms, stadiums, parking lots, and drive-through points.



\begin{figure}[!h]
\centering
\includegraphics[width=10cm]{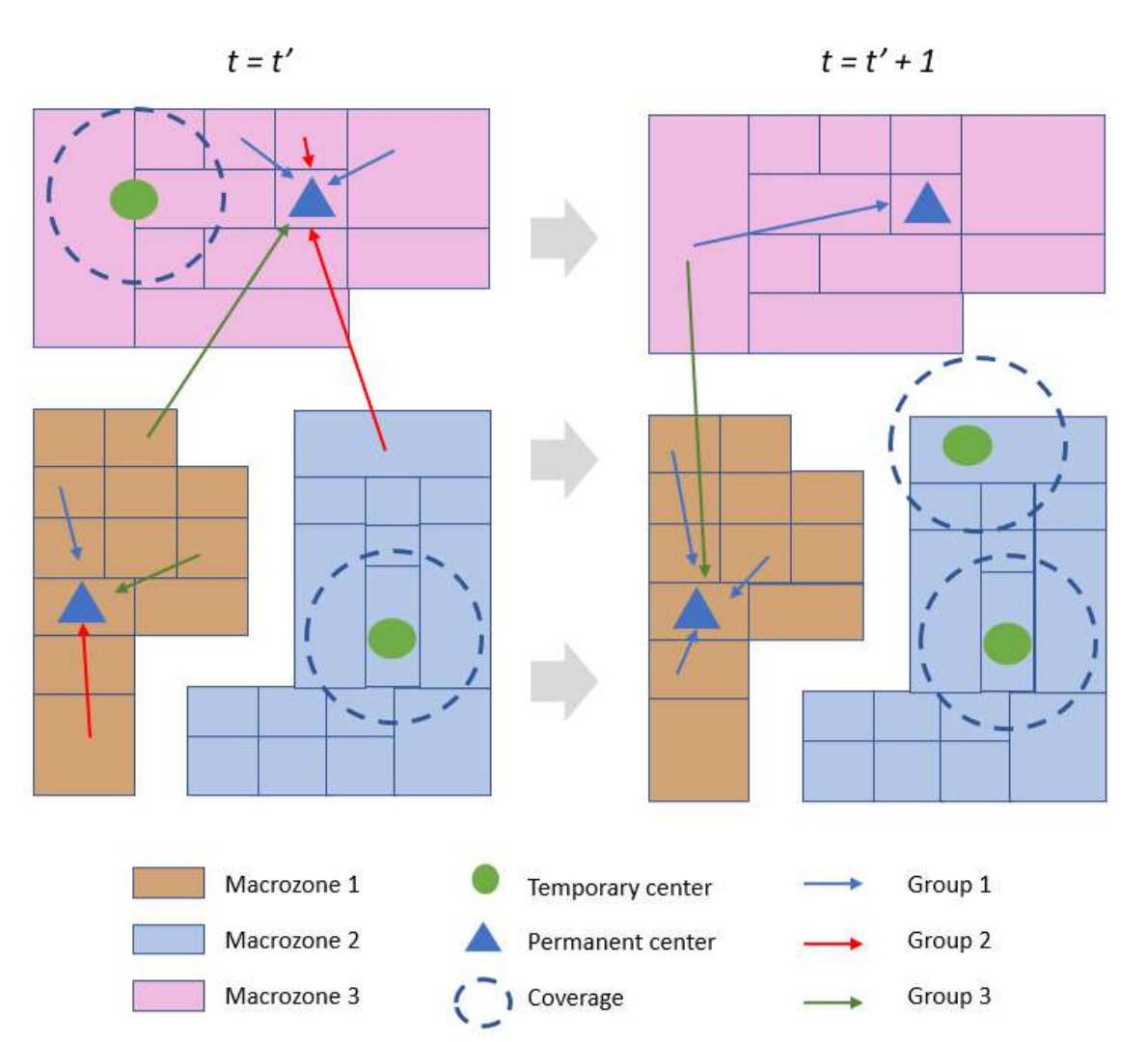}
\caption{\footnotesize A feasible solution for two days $t=t'$ and $t=t'+1$ is displayed. Rectangles represent neighborhoods, and their color indicates the zone in the city (or macrozone). Triangles denote permanent vaccination centers. Circles indicate temporary vaccination centers. Colored arrows show the target groups. The segmented-blue circles show the coverage area of each temporary center (adjacent zones).  \label{figure2} }
\end{figure}

Figure \ref{figure2} presents an area where our problem can be illustrated. We consider a region in which we want to vaccinate the population. This area is divided into macrozones. The example contains a set $K$ of three macrozones, which are denoted by different colors. Moreover, each macrozone $k \in K$ has a set of neighborhoods $L_k$, which is disjoint from the sets from the other macrozones. It also shows three macrozones with a different number of neighborhoods; e.g., the macrozone on top (color pink) has thirteen (rectangles) neighborhoods. We assume that the effective demand, which is considered at the neighborhood level, is known. The demand must be satisfied completely, and it is divided into a set of priority groups, e.g., ill and disabled people, older people, children, and adults. Figure \ref{figure2} shows three priority groups, marked with a different arrow color. The arrow direction indicates where each group is served. The permanent centers can provide vaccines to people from any place in the geographical region.

The problem assumes that priorities exist in those groups. People can be vaccinated using either permanent centers (triangles) or temporary centers (circles). It means that the coordination of these resources during the vaccination operations is crucial. Moreover, to execute the vaccination process, a planning horizon is considered. In Figure  \ref{figure2} a feasible solution for two consecutive days ($t = t'$ and $t = t' + 1$) is presented.  Each day additional staff to those in the permanent centers work in a temporary center located at a site selected from a possible set of locations. \cite{Lee2020PoliciesKorea} collected evidence from rapid full testing campaigns using low-cost infrastructure into the drive-through and walk-through systems, showing some advantages of using temporary centers.

\subsection{Mathematical formulation for the MMV problem}
\label{mathematical_formulation}

The mathematical formulation for the problem has two objectives. The first focused  on executing the vaccination plan as soon as possible while prioritizing high-risk groups (similar to minimizing the makespan in a scheduling context). The second defines an economic objective associated with using temporary centers, representing an additional cost for the vaccination campaign. This cost considers payments for contracting additional staff. Regarding the constraints, each feasible solution is defined by the next three sets of constraints:  the first one balances the flow of vaccines among the levels (distribution center / macrozone / neighborhood). The second one is focused on capacity. The third one is related to temporary centers' locations. We call our problem Multi-Modal Vaccination (MMV) problem, which has as a particular case the set covering problem,  a traditional NP-Hard problem. Therefore, MMV is also NP-Hard.

The model considers the following sets: $K$ is the set of macrozones of the region under study. $L_{k}$ is the set of neighborhoods from macrozone $k \in K$. Note that each neighborhood can belong to one macrozone, this means that $L_{1}, L_{2}, .. L_{k}, .., L_{|K|}$ are disjoint sets. $T$ is the set of days in which to carry out the vaccination campaign. $P$ is the set of population groups. $I$ is the set of permanent centers for vaccination. $J$ is the set of temporary centers for vaccination. $N_{l}$ is the set of neighborhoods that can be served by a temporary center from the neighborhood $l \in \left\lbrace L_{k}: k \in K \right\rbrace $.

\noindent {The optimization model includes the following parameters:} \\ [0.5em]
\begin{tabular}{l p{10cm}}
    $pop_{lp}$ & Demand for vaccination of the group $ p \in P $ in neighborhood $ l \in L_{k} $ of macrozone $ k \in K $.  \\  $C_i$ & Vaccination capacity per day of permanent center $ i \in I $. \\ 
    $D_j$ & Vaccination capacity per day of temporary center $ j \in J $ to serve each macrozone. \\ 
    $A_{t}$ & Maximum supply from the vaccine distribution center for the day $ t \in T $.\\
    $r_{p}$ & Index of risk level of group $ p \in P $.  \\ 
    $\epsilon_{p}$ & Daily variation in the risk level of people belonging to group $ p \in P $.  \\ 
    $mc$ & Cost of using a temporary center (staff, facilities, etc.) 
\end{tabular}

\noindent {The decision variables are: } \\ [0.5em]
\begin{tabular}{l p{10cm}}
    $\phi_{lpti}$ & Number of people in group $ p \in P $ in neighborhood $ l \in \left\lbrace L_{k}: k \in K \right\rbrace $, served on day $ t \in T $ by permanent center $ i \in I $. \\
    $\gamma_{lptj}$ & Number of people in group $ p \in P $ in neighborhood $ l \in \left\lbrace L_{k}: k \in K \right\rbrace $, served on day $ t \in T $ in temporary center $ j \in J $. \\
    $y_{jtl}$ & Binary variable that turns on when temporary center $ j \in J $ is installed in neighborhood $ l \in L_{k} $ of macrozone $ k \in K $ on day $ t \in T $.\\
    $v_{jptl}$ & Binary variable that turns on when temporary center $ j \in J $ serves population group $ p \in P $ in neighborhood $ l \in L_{k} $ of macrozone $ k \in K $ on day $ t \in T $. \\
\end{tabular}


Based on this notation, a bi-objective mixed-integer linear programming model for the MMV problem is formulated as follows: 

\begin{align}
  f=& (f_1, f_2) \label{fos} \\ 
 \min f_1 = &\sum_{p \in P, t \in T}  \left[  (1-r_{p})(1+\epsilon_p)^t  (\sum_{l \in \left\lbrace L_{k}: k \in K \right\rbrace, i \in I} \phi_{lpti} + \sum_{l \in \left\lbrace L_{k}: k \in K \right\rbrace, j \in J} \gamma_{lptj}) \right] \label{f_1} \\
 \min f_2 =& \sum_{j \in J, t \in T ,l \in \left\lbrace L_{k}: k \in K \right\rbrace} mc\cdot y_{jtl} \label{f_2}
 \end{align}
 
 \noindent subject to: 
  \begin{align}
   \
    \sum_{l \in \left\lbrace L_{k}: k \in K \right\rbrace, i \in I} \phi_{lpti} + \sum_{l \in \left\lbrace L_{k}: k \in K \right\rbrace, j \in J} \gamma_{lptj}&\leq A_t && \forall t \in T \label{r4} \\
    \sum_{t \in T, i \in I} \phi_{lpti} + \sum_{t \in T, j \in J} \gamma_{lptj} &\geq  pop_{lp} && \forall l \in \left\lbrace L_{k}: k \in K \right\rbrace, \thinspace p \in P \label{r5} \\
    \sum_{l \in \left\lbrace L_{k}: k \in K \right\rbrace, p \in P}\gamma_{lptj} &\leq D_{j} && \forall j \in J, t \in T \label{r6}\\
    \sum_{l \in \left\lbrace L_{k}: k \in K \right\rbrace, p \in P}\phi_{lpti} &\leq C_{i} && \forall i \in I, t \in T \label{r7}\\
    \sum_{ r \in N_{l}} y_{jtr} - v_{jptl} &\geq 0 &&  \forall p \in P, \thinspace j \in J, t \in T, \thinspace l \in \left\lbrace L_{k}: k \in K \right\rbrace \label{r8} \\
    v_{jptl} \cdot pop_{lp} - \gamma_{lptj}  &\geq 0 && \forall p \in P, \thinspace j \in J, \thinspace l \in \left\lbrace L_{k}: k \in K \right\rbrace \thinspace ,t \in T \label{r9}\\
    \sum_{l \in \left\lbrace L_{k}: k \in K \right\rbrace} y_{jtl} & \leq 1 && \forall j \in J, t \in T   \label{r10} \\
     \phi_{lpti}, \thinspace \gamma_{lptj} & \geq 0 && \forall l \in \left\lbrace L_{k}: k \in K \right\rbrace, p \in P, t \in T, i \in I, j \in J \label{r11} \\
     y_{jtl}, v_{jptl} &\in \{ 0, 1 \}  && \forall j \in J, p \in P, t \in T, l \in \left\lbrace L_{k}: k \in K \right\rbrace \label{r12} 
\end{align}

 Regarding objective functions, $f_1$ ensures that the high-risk population receives its vaccines as soon as possible. $f_1$ is a non-normalized weighted sum of vaccinated people whose weights are modeled by $(1-r_{p})(1+\epsilon_p)^t$, which also can be understood as a temporal prioritization function. This function has an exponential growth on time ($t$), where its growth factor $(1+\epsilon_p)$ and the initial prioritization value $(1-r_{p})$ depend on the group $p$. Due to $f_1$ is minimized, low values for the weights are prioritized for each person (i.e., high values for the risk $r_p$ and small values for $t$). Hence, when higher values of $\epsilon_p$ are considered for a group, the model finds to dispense vaccines as soon as possible, avoiding the exponential prioritization function growth in this group. $\epsilon_p$ is a parameter based on expert knowledge that quantifies the daily variation in the initial risk for each group. Let's see an example, a person $p_1$ belonging to a high-risk group can have a weight of $(1-0.8)(1+0.3)^1 \approx 0.26 $ at day 1 and $(1-0.8)(1+0.3)^5 \approx 0.75$ at day 5, and other person $p_2$ of middle risk can have a weight of $(1-0.2)(1+0.1)^1  \approx 0.88$ on the first day and $(1-0.2)(1+0.1)^5 \approx 1.23$ on the fifth day. The decision here is to select which person is vaccinated each day, assuming two days (first day or fifth day) and one vaccine per day. Then, the more convenient action is to vaccine $p_1$ on the first day and $p_2$ on the fifth day, scoring 1.49; the other alternative scores 1.65. This is the prioritization incorporated into $f_1$. On the other hand, as previously mentioned, temporary centers offer an opportunity to improve service quality by minimizing congestion or speeding up the campaign. However, these resources represent a relevant cost (monetary cost or opportunity cost) for the health network; therefore, it is also essential to consider the efficient use of those resources, which is modeled by (\ref{f_2}). 

Constraints (\ref{r4}) limit the people vaccinated to the maximum number of available vaccines at each period. Constraints (\ref{r5}) establish that all people from the neighborhoods should be vaccinated at the campaign's end. Constraints (\ref{r6}) claimed that the available capacity of vaccines is respected at each period and at each temporary center. Similar constraints (\ref{r7})are defined for the capacity of vaccines at the permanent centers. Constraints (\ref{r8}) claim that only people from neighborhoods in the coverage area of temporary centers are attended. Constraints (\ref{r9}) establish that the number of people vaccinated in temporary centers must be lower than or equal to the demand of each group in every neighborhood. Constraints (\ref{r10}) restrict each temporary center to be opened each day $t$ in an unique neighborhood. Finally, constraints (\ref{r11}) and (\ref{r12}) specify the nature of the variables.

To solve this bi-objective problem, we apply a traditional weighted sum (blended) approach, then, the weighted sum objective is: 

\[ \min Z = \alpha f_1^\ast + (1-\alpha) f_2^\ast, \]

\noindent  $ 0\leq \alpha \leq 1$, and $f_1^\ast$, $f_2^\ast$ are normalized objective functions. 
This approach has the important property that the weighted sum $Z$ problem has the same computational complexity and needs the same computational effort to solve as the single objective version of the multi-objective MVV problem. However, an important issue in the application of this method is what values of $\alpha$ must be used. The
values depends on the importance of each objective in the particular problem and a scaling factor. The scaling effect can be avoided somewhat by normalizing the objective functions. \cite{Marler2010TheInsights} and \cite{Deb2014Multi-objectiveOptimization} discuss several practical issues about the implementation of this method. To normalize the objectives we use the standard 0-1 min-max scaling that is defined as follows:

\[ f^\ast = \dfrac{f - f_{\min}}{f_{\max} - f_{\min}}, \] 

\noindent where $f_{\min}$ and $f_{\max}$ are obtained from single-objective optimization models\footnote{\textbf{Case $\min f$:} Let $f_{\min}$ be the optimal solution in the single-objective problem. Let $f_{\max}$ be the unbounded constraint value from an optimization model considering the other objective. \textbf{Case $\max f$}. Let $f_{\max}$ be the optimal solution in the single-objective problem. Let $f_{\min}$ be the unbounded constraint value from an optimization model considering the other objective.}.

\subsection{Illustrative example}
\label{3.3}
To illustrate some properties of the previously presented model, we use the following instance. We consider a region that is divided into four macrozones, $K= \left\lbrace 1,2,3,4\right\rbrace $, using a set of five neighborhoods for each one. Moreover, the demand is grouped into three population segments, $P= \left\lbrace A, B, C\right\rbrace $. The risk level for those groups is defined as follow: $A$, high risk ($r_A= $0.8); $B$, medium risk ($r_B=$ 0.5); and $C$, low risk ($r_C =$ 0.28). The $\epsilon_p, p \in P$, values are $\epsilon_A=$ 0.06, $\epsilon_B=$0.025, and $\epsilon_C=$ 0.015. Hence, the temporal prioritization function for each group can be expressed as follows:
\begin{align}
    \text{group A} = (1-r_{A})(1+\epsilon_A)^t = (0.20)(1.060)^t, \\
    \text{group B} = (1-r_{B})(1+\epsilon_B)^t = (0.50)(1.025)^t,  \\
    \text{ and }  \text{group C} = (1-r_{C})(1+\epsilon_C)^t = (0.72)(1.015)^t.
\end{align}

\begin{figure}
\centering
\includegraphics[width=12cm]{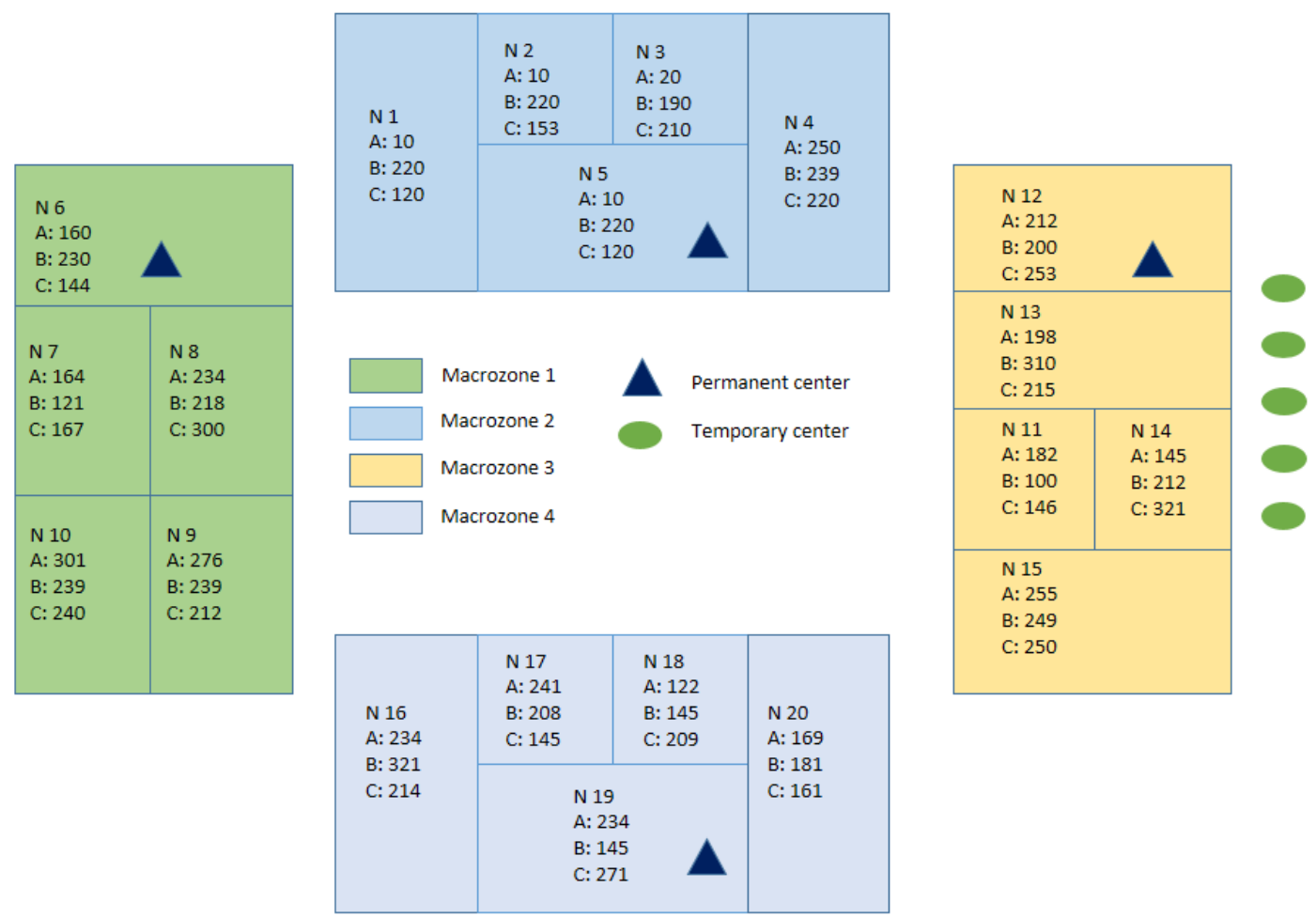}
\caption{\footnotesize Geographic adjacency of neighborhoods and their population for each priority group.}
\label{figure33}
\end{figure}

The aggregated demand for each of those groups are $A= 3,657, B=3,906$ and $C= 4,051$, dis-aggregated requirement for each neighborhoods and macrozones are presented in Figure \ref{figure33}. Regarding temporary centers, we consider five centers, their capacity is 37 vaccines, and they can cover adjacent neighborhoods. The cost per day of use is \$350; it represents a variable transport cost and an opportunity cost by reducing the health center's capacity from where they come. Our example considers a permanent center for each macrozone; its capacity is 150 vaccines per day. Lastly, a 20-day horizon is considered for serving this demand.

When looking at Table \ref{tablepareto}, we can observe the trade-off between objective functions. For example, if we prioritize $f_1$ (higher values for $\alpha$), then the percentage of people served by temporary centers increases (see column $\mathcal{P}$), and the vaccination period is shorter ($\mathcal{D}$ down to 16 days). On the other hand, if we prioritize $f_2$ (lower values for $\alpha$), then the cost of temporary centers is more important, the percentage of people vaccinated at temporary centers is smaller, and the vaccination period is more extended ($\mathcal{D}$ up to 20 days). Specifically, when $\alpha$ assumes values 0 and 0.4, it is possible to have a slight improvement in the use of the resources; that is, the value of $\mathcal{P}$ is slightly increased,  maintaining the lowest cost of temporary centers. The $\mathcal{P}$ value increased from 5.62\% to 5.78\%. The run times can be seen in the last column of Table \ref{tablepareto}. The model is easier to solve for the mono-objective optimization problems ($\alpha \in \{0,1\}$).  
The computer used to run the models has a processor AMD Rayzer 5 3450U with Radeon Vega Mobile Gfx 2.10 GHz. and 8 GB of RAM. Based on previous experiments, we tuned the Gurobi's \textit{MIRCuts} parameter for an aggressive cuts generation and the Branch \& Cut algorithm was set to prioritize optimality over feasibility. The code and data will be made available at \url{https://github.com/fmd0061566/biobjectivevaccines}.

\begin{table}[]
\centering
\footnotesize
\begin{tabular}{|c|cccc|c|c|c|}
\hline
\multirow{2}{*}{$\alpha$} & \multicolumn{4}{c|}{Objective Functions Values}                                                     & \multirow{2}{*}{{$\mathcal{P}$}} & \multirow{2}{*}{{$\mathcal{D}$}} & \multirow{2}{*}{Runtime (s)} \\ \cline{2-5}
                       & \multicolumn{1}{c|}{$f_1^*$}   & \multicolumn{1}{c|}{$f_1$}       & \multicolumn{1}{c|}{$f_2^*$}   & $f_2$        &                    &                    &                              \\ \hline
0                      & \multicolumn{1}{c|}{0.902} & \multicolumn{1}{c|}{7,605.567} & \multicolumn{1}{c|}{0.000} & 4,900.000  & 5.62\%             & 20                 & 119                          \\
0.4                    & \multicolumn{1}{c|}{0.423} & \multicolumn{1}{c|}{7,290.170} & \multicolumn{1}{c|}{0.002} & 4,955.732  & 5.78\%             & 20                 & 600                          \\
0.6                    & \multicolumn{1}{c|}{0.132} & \multicolumn{1}{c|}{7,099.130} & \multicolumn{1}{c|}{0.306} & 14,116.162 & 15.86\%            & 18                 & 600                          \\
0.8                    & \multicolumn{1}{c|}{0.011} & \multicolumn{1}{c|}{7,019.060} & \multicolumn{1}{c|}{0.588} & 22,595.792 & 24.78\%            & 16                 & 600                          \\
1                      & \multicolumn{1}{c|}{0.000} & \multicolumn{1}{c|}{7,012.000} & \multicolumn{1}{c|}{1.000} & 35,000.000  & 28.60\%            & 16                 & 1                            \\ \hline
\end{tabular}
\caption{\footnotesize Model results for different $\alpha$ values. $\mathcal{P}$ represents the percentage of patients vaccinated by temporary centers. $\mathcal{D}$ is the number of days to finish the campaign. The runtime for each $\alpha$ value, is also shown. }
\label{tablepareto}
\end{table}

 For three values of $\alpha$, their vaccination plan advance are presented in Figure \ref{ilustrativo}. It shows the vaccination plan for risk groups. For $\alpha=0$, Figure \ref{ilustrativo}-(a), it can be seen that there is no prioritization for risk groups. Note that, every day, people from each group are vaccinated. On the contrary, for $\alpha=0.8$, note that there is a clear prioritization for risk groups. Risk group A finishes the vaccination on day 6, risk group B finishes the vaccination on day 11, and day 18 for risk group C. This effect is intensified for $\alpha=1$, noting that the campaign ended on day 16 in this case, and on day 18 for $\alpha=0.8$.

\begin{figure}
\centering
\subfloat[]{\includegraphics[width=9.5cm]{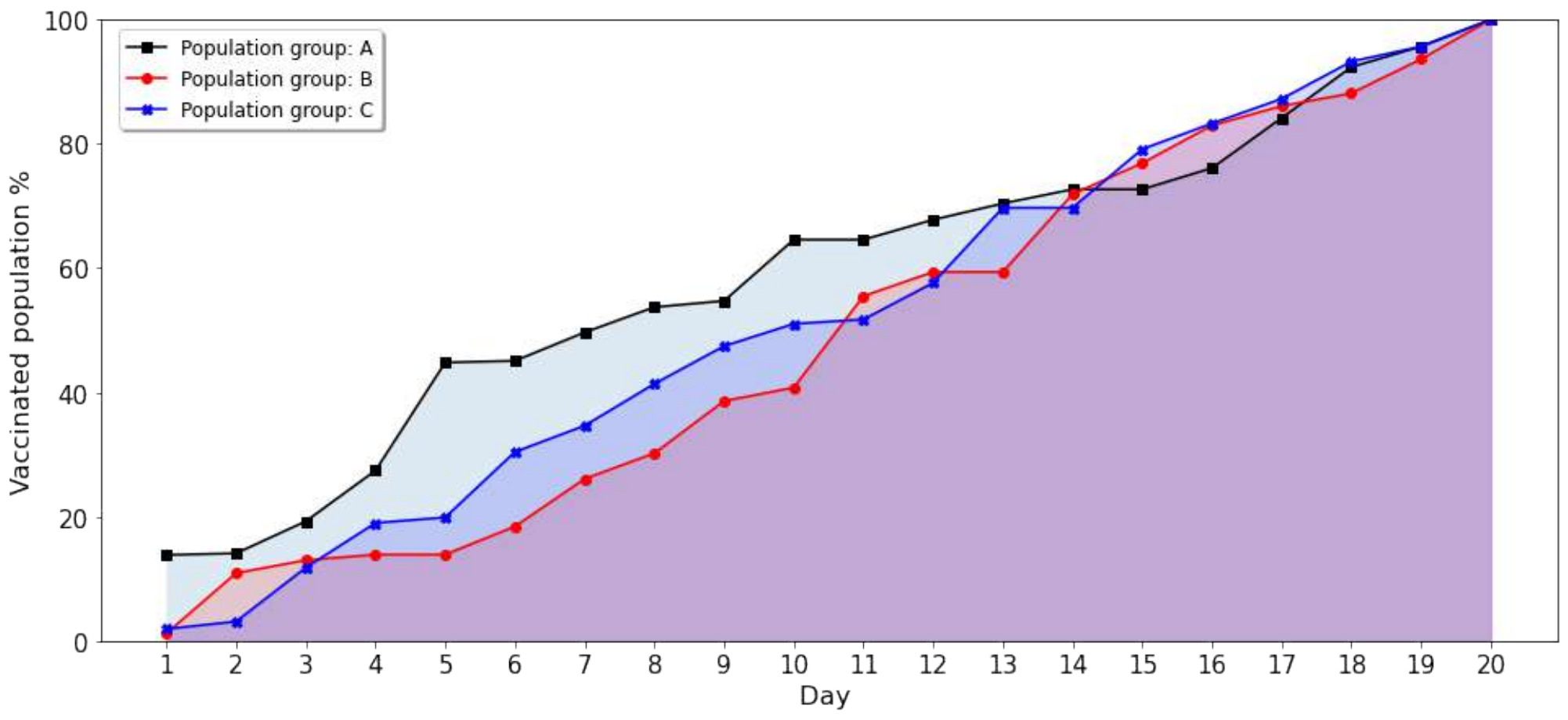}}\hspace{1mm}
\subfloat[]{\includegraphics[width=9.5cm]{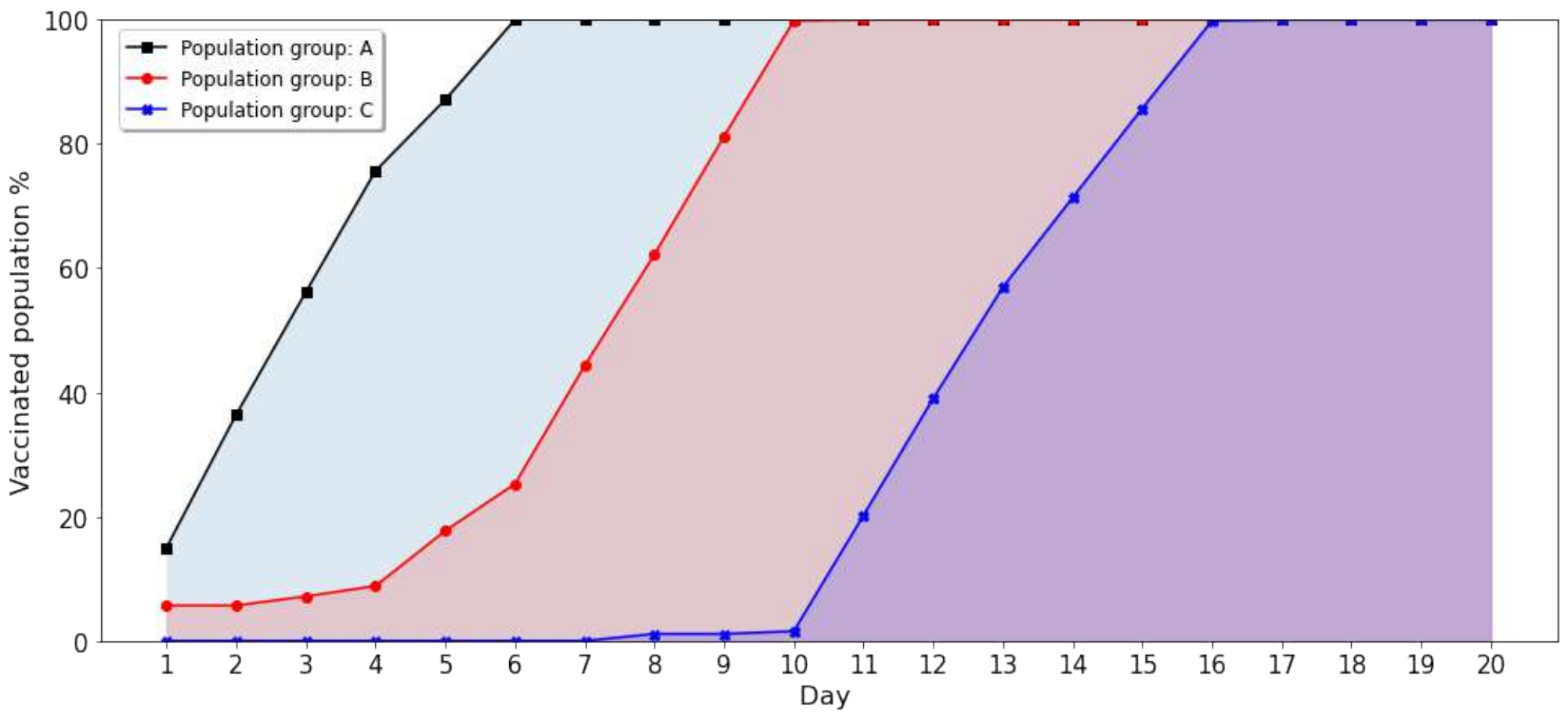}}\hspace{1mm}
\subfloat[]{\includegraphics[width=9.5cm]{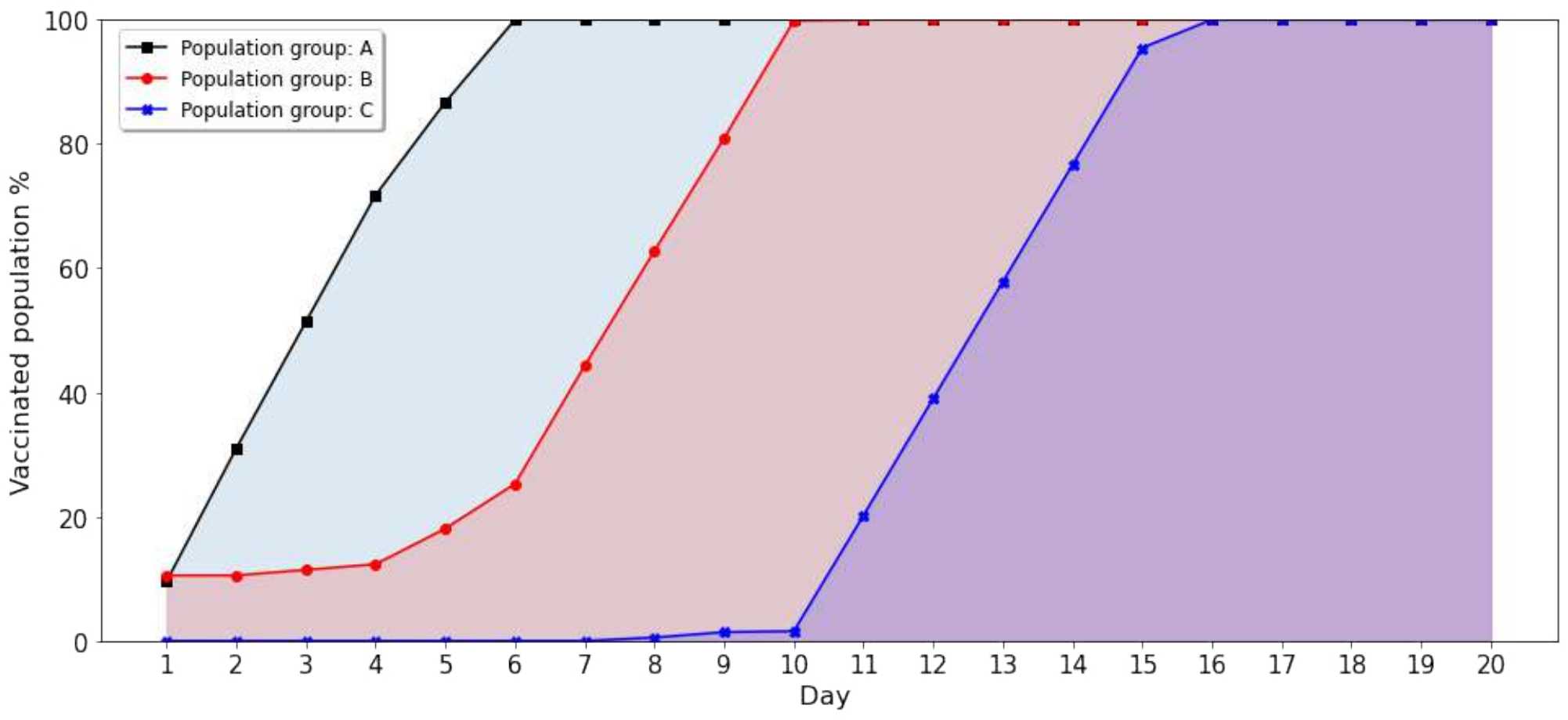}}\hspace{1mm}
\caption{{ \footnotesize Percentage of vaccinated people on day $t$ by group. When $\alpha=$ 0 (a), $\alpha=$ 0.8 (b), and $\alpha=$ 1 (c).}\label{ilustrativo}}
\label{ilustrativo}
\end{figure}

\begin{table}
\centering
\footnotesize
\begin{tabular}{c|ccccccccccc|ccc|c}
\cline{2-15}
                          & \multicolumn{11}{c|}{Neighborhoods, $j=1$}             & \multicolumn{3}{c|}{Groups} &                            \\ \hline
\multicolumn{1}{|c|}{Day} & 3  & 6  & 8  & 9  & 11 & 12 & 13 & 14 & 15 & 16 & 19 & 1       & 2       & 3       & \multicolumn{1}{c|}{Total} \\ \hline
\multicolumn{1}{|c|}{1}   &    & 37 &    &    &    &    &    &    &    &    &    & 37      &         &         & \multicolumn{1}{c|}{37}    \\
\multicolumn{1}{|c|}{2}   &    &    & 10 & 27 &    &    &    &    &    &    &    & 37      &         &         & \multicolumn{1}{c|}{37}    \\
\multicolumn{1}{|c|}{3}   &    & 37 &    &    &    &    &    &    &    &    &    &         & 37      &         & \multicolumn{1}{c|}{37}    \\
\multicolumn{1}{|c|}{4}   &    &    &    &    &    &    &    &    &    &    & 37 & 37      &         &         & \multicolumn{1}{c|}{37}    \\
\multicolumn{1}{|c|}{5}   & 37 &    &    &    &    &    &    &    &    &    &    &         & 37      &         & \multicolumn{1}{c|}{37}    \\
\multicolumn{1}{|c|}{6}   &    &    &    &    &    &    &    &    &    &    & 37 & 37      &         &         & \multicolumn{1}{c|}{37}    \\
\multicolumn{1}{|c|}{7}   &    &    &    &    &    &    &    &    &    & 37 &    &         & 37      &         & \multicolumn{1}{c|}{37}    \\
\multicolumn{1}{|c|}{8}   &    &    &    &    &    &    & 37 &    &    &    &    &         & 37      &         & \multicolumn{1}{c|}{37}    \\
\multicolumn{1}{|c|}{9}   &    &    &    &    & 7  & 30 &    &    &    &    &    &         & 37      &         & \multicolumn{1}{c|}{37}    \\
\multicolumn{1}{|c|}{10}  &    & 37 &    &    &    &    &    &    &    &    &    &         & 37      &         & \multicolumn{1}{c|}{37}    \\
\multicolumn{1}{|c|}{11}  &    &    &    &    &    &    &    & 37 &    &    &    &         &         & 37      & \multicolumn{1}{c|}{37}    \\
\multicolumn{1}{|c|}{12}  &    &    &    &    &    &    &    &    & 37 &    &    &         &         & 37      & \multicolumn{1}{c|}{37}    \\
\multicolumn{1}{|c|}{13}  &    &    &    &    &    &    &    &    &    &    & 37 &         &         & 37      & \multicolumn{1}{c|}{37}    \\ \hline
\end{tabular}
\caption{Number of vaccines delivered daily for different neighborhoods through the temporary center $j=1$. with $\alpha=$ 0.8}
\label{tabla_ilustrativo}
\end{table}

To complement the analysis of the illustrative example, we present Table \ref{tabla_ilustrativo}. It shows the deployment of a temporary vaccination center throughout the campaign. First of all, it should be noted that the full capacity of this center is used, 37 vaccinations per day. On the other hand, it can be observed that the most at-risk groups (groups A and B) are prioritized during the first days of the vaccination campaign, leaving the group with the lowest priority (group C) for the campaign's last days. Finally, it can be seen that once installed in a sector for a given day, it can serve more than one neighborhood according to service level criteria. Table \ref{tabla_ilustrativo} shows this effect for days 2 and 9. To recall the geographic distribution, see Figure \ref{figure33}. 

In this illustrative example, we can observe the variety of solutions provided by our bi-objective model. The planner can decide the campaign's prioritization, considering the availability of resources and the level of risk in the environment. Sections \ref{section4} and \ref{section5} present a more in-depth study of the model solutions using a real data set.

\section{A case study: Vaccination campaign planning in San Bernardo }\label{section4}

The case study's goal is to show how the model help to plan a single-dose vaccination campaign. The model defined in Section 3 is illustrated in Chile; specifically, San Bernardo is located in the South West of the Metropolitan Region. San Bernardo, with an estimated population of 334,836 inhabitants \cite{GajardoPolanco2019Region2015-2035}, is the sixth most populated commune in the Metropolitan Region and represented 23 percent of the total population of the South Metropolitan Health Service in 2019.

The case is developed using current population data from the commune of San Bernardo in Chile and the number of neighborhoods and permanent vaccination centers. The case study uses information from the "Pre-census 2017" by the National Institute of Statistics \cite{DepartamentodeGeografia2019DivisionSantiago}. The macrozones, set $K$, were defined as the geographic grouping units. Then, inside them, the neighborhoods, set $L_k$, were identified. Table  \ref{tab:c12} shows the macrozones' total urban and rural population of San Bernardo. The number of Primary health centers corresponds to the number of permanent vaccinations centers (model set $I$).

\begin{figure}
\centering
\includegraphics[width=12cm]{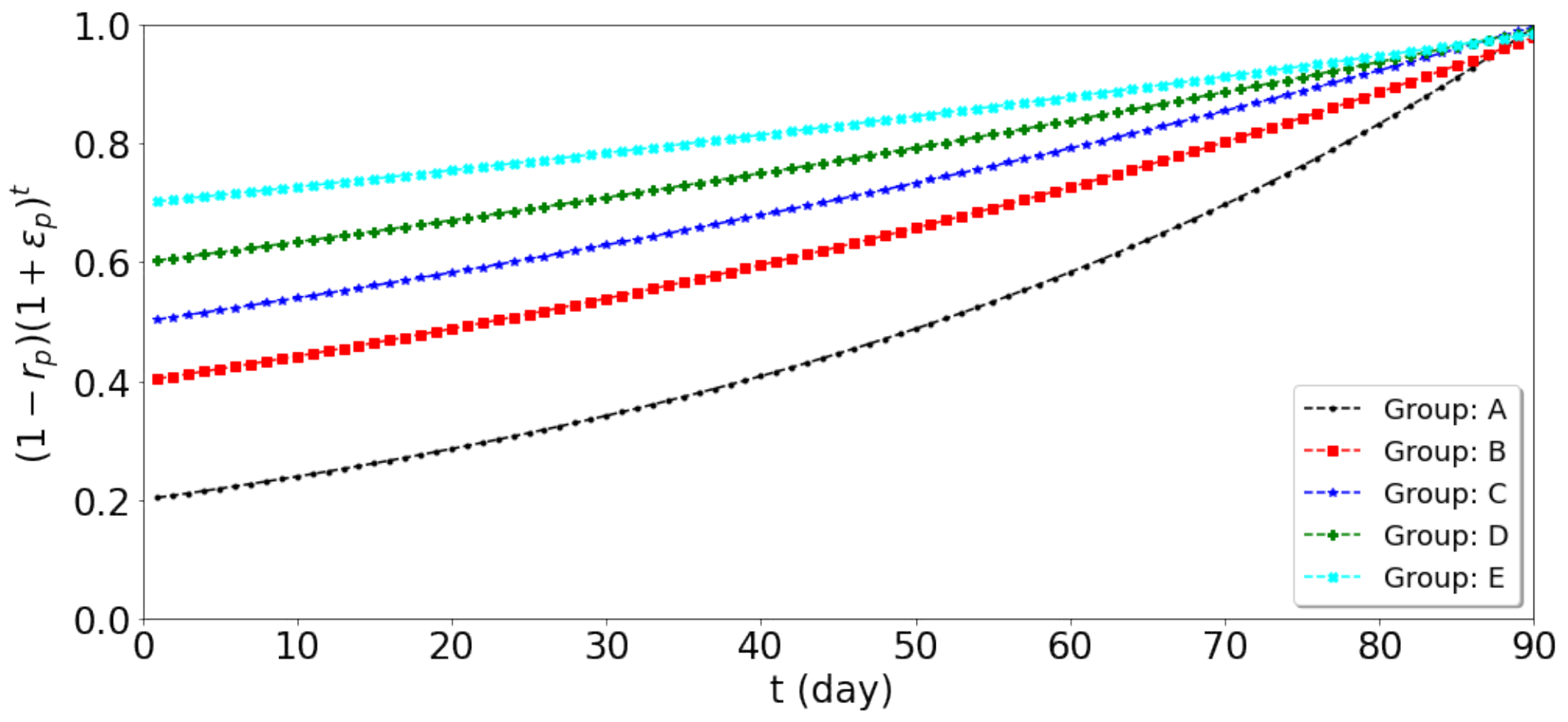}
\caption{\footnotesize Behavior of the temporal prioritization function for each group, throughout the vaccination campaign.   }
\label{f:r}
\end{figure}

The case study considers five target groups. These groups are characterized as follows: group A, people aged 65 years and over; group B, those with chronic illnesses from 11 to 64 years; group C: pregnant women; group D, children from 1st to 5th grade; and group E: children from 6 months to 5 years. In Figure \ref{f:r} we present the previously explained (Section 3) priority function for each of these groups is shown. Note that the prioritization function for a higher risk person (e.g., persons in group A) has a steeper exponential growth concerning a person from less risky groups (e.g., persons in group E). Therefore, since the objective function is minimizing, it is more attractive to vaccinate persons in group A in the first days of the vaccination campaign compared to persons in group E.


\begin{figure}
\centering
\includegraphics[width=13cm]{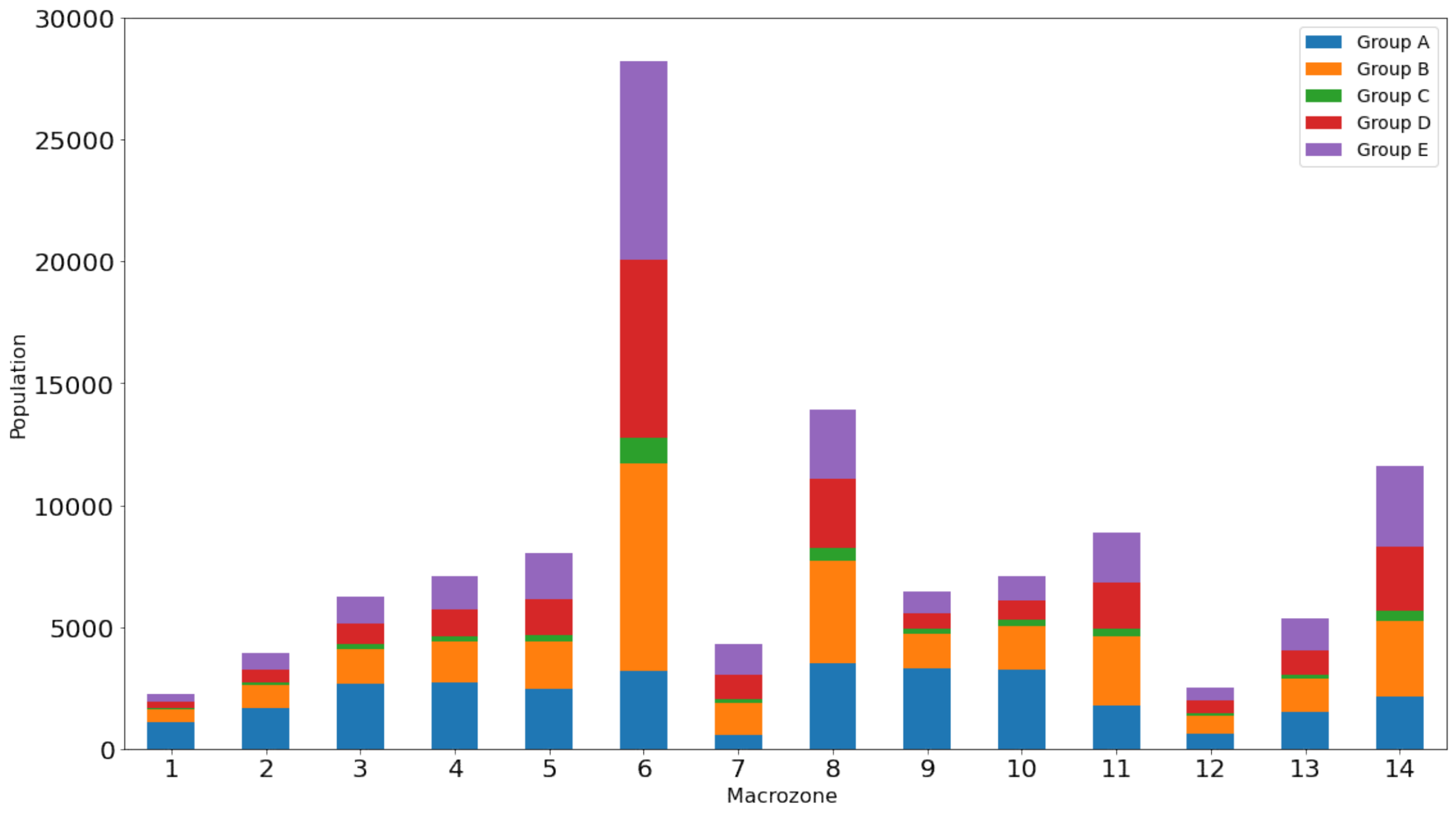}
\caption{Distribution of population groups in macrozones.  Source: \cite{DepartamentodeGeografia2019DivisionSantiago}
\label{fy}}
\end{figure}

Lastly, the parameter $pop_{lp}$, which represents the demand for vaccines, was estimated based on data from the "Pre-census 2017" from the National Institute of Statistics \cite{DepartamentodeGeografia2019DivisionSantiago} estimation of the population groups within each of the neighborhoods and macrozones. This information is presented in Figure \ref{fy}. 


\begin{table}[]
\centering
\footnotesize
\begin{tabular}{|c|c|c|c|c}
\hline
\multicolumn{5}{|c|}{San Bernardo}                                \\ \hline
id & Macrozone & Surface Km$^2$ & Population Census 2017 & \multicolumn{1}{c|}{Households Census 2017} \\ \hline
1 & O'Higgins             & 1.2  & 5,100  & \multicolumn{1}{c|}{2,045}  \\
2 & Escuela de Infantería & 0.9  & 9,204  & \multicolumn{1}{c|}{2,664}  \\
3 & Calderón de La Barca  & 1.0    & 14,182 & \multicolumn{1}{c|}{4,036}  \\
4 & Santa Marta           & 1.2  & 16,403 & \multicolumn{1}{c|}{4,390}  \\
5 & Hospital              & 1.4  & 18,994 & \multicolumn{1}{c|}{5,634}  \\
6 & Cerro Negro           & 19.7 & 79,725 & \multicolumn{1}{c|}{23,463} \\
7 & Maestranza            & 1.7  & 11,681 & \multicolumn{1}{c|}{4,097}  \\
8 & Nos                   & 21.7 & 38,492 & \multicolumn{1}{c|}{11,404} \\
9 & Nogales               & 2.2  & 14,387 & \multicolumn{1}{c|}{4,737}  \\
10 & Tejas de Chena        & 1.8  & 17,591 & \multicolumn{1}{c|}{5,406}  \\
11 & Chena                 & 34.8 & 21,150 & \multicolumn{1}{c|}{7,381}  \\
12 & Lo Herrera            & 51.1 & 7,051  & \multicolumn{1}{c|}{2,135}  \\
13 & Estación              & 0.9  & 12,915 & \multicolumn{1}{c|}{3,705}  \\
14 & Los Morros            & 13.2 & 29,730 & \multicolumn{1}{c|}{8,951}  \\
-& Stragglers            & -    & 698   & \multicolumn{1}{c|}{213}   \\ \hline
\end{tabular}
\caption{Detail of the population and households of the commune of San Bernardo. Source: \cite{DepartamentodeGeografia2019DivisionSantiago}}
\label{tab:c12}
\end{table}

To give a more comprehensive vision of the vaccination campaign in San Bernardo, two scenarios were defined. We seek to assess the performance of our model when there is a low/high (Scenario 1/Scenario 2) capacity in permanent centers. The capacities summary is presented in Table \ref{tab:capacidades} and the details of each scenario are presented below. 

\begin{itemize}
    \item Scenario 1 ($s_1$): Five temporary vaccination centers are available per day, with a daily capacity of 200 vaccines for each one ($D_j$). The daily cost of a temporary center is \$350 ($mc$). This scenario considers a supply capacity ($A_t$) equal to 1,800 vaccines. The total capacity of permanent centers in $s_1$ is 623 vaccines. Figure \ref{f:r} presents the temporal prioritization function for each group.
    \item Scenario 2 ($s_2$): The vaccination capacity ($D_j$) and cost of temporary centers ($mc$) remain the same as in Scenario 1. The changes are in both the supply capacity ($A_t$) and the permanent vaccination centers capacity ($C_i$).
\end{itemize}

\begin{table}
\centering
\footnotesize
\begin{tabular}{c|c|c|}
\cline{2-3}
                         & $s_1$   & $s_2$   \\ \hline
\multicolumn{1}{|c|}{$A_t$} & 1,800 & 2,000 \\ \hline
\multicolumn{1}{|c|}{$\sum_{i \in I}C_i$} & 623  & 1,000 \\ \hline
\multicolumn{1}{|c|}{$\sum_{j \in J}D_j$} & 1,000 & 1,000 \\ \hline
\multicolumn{1}{|c|}{$\sum_{l \in L_k,p \in P}pop_{pl}$} & 115,800  & 115,800  \\\hline
\multicolumn{1}{|c|}{$mc$} & 350  & 350  \\ 
\hline

\end{tabular}

\caption{\footnotesize Daily supplied vaccines, total vaccination capacities by center type, total demand and cost of temporary centers.}
\label{tab:capacidades}
\end{table}

\section{Results and managerial insights}
\label{section5}

This section discusses several viewpoints about the vaccination campaign's solutions provided by our proposed model. We test this model using the two scenarios presented above ($s_1$ and $s_2$). This section is organized as follows. Section \ref{5.1} gives the implementation setting. Section \ref{5.2} presents a comparison between the two scenarios and provides the $\alpha$-value impact on two key performance indicators (prioritization of groups and duration). In Section \ref{5.3}, a numerical discussion about objective functions and different values of $\alpha$ is presented. In Section \ref{5.4}, experiments and results about the vaccination campaign's deployment in different macrozones are shown. In Section \ref{5.5}, the level of risk and temporary centers use are contrasted. Then, in Section \ref{5.6}, to minimize the congestion of higher priority groups, a new constraint is added, and the results from this variant are presented. Section \ref{5.7} shows an extension of the proposed model to deal with uncertainty in vaccine demand. Lastly, Section \ref{5.8} gives a set of useful managerial insights from our research.




\subsection{Implementation and code}
\label{5.1}
The following experiments were implemented using the Python 3.6 language and Gurobi 9.0.2 as a solver. The experiments were run on a computer with an Intel Core i7-8700 CPU with a base operating frequency of 3.20 GHz. 64 GB of RAM, and twelve processors. In the following experiments, we use the same solver setting as in the Illustrative Example.

\subsection{Objectives' prioritization impact on the vaccination campaign}
\label{5.2}
In this section, different analyses are carried out to demonstrate the main advantages of using this tool for decision-making. The two scenarios presented in the previous section are considered. Taking advantage of our model's bi-objective approach, after a preliminary experimentation, two $\alpha$ values are evaluated to study the prioritization of one objective over another in order to appreciate its impact on the vaccination campaign. 


For each scenario, two vaccination schedules are presented in Figure \ref{e1}.  For $s_1$, we show the vaccination for $\alpha=0.2$ in Figure \ref{e1}-a and for $\alpha=0.98$ in Figure \ref{e1}-b. In the former, the accumulated percentage of vaccinated people by group is shown (e.g., group A is completely vaccinated by day 44 of the campaign) where we can see that the population groups' attention is prioritized. Also, when there is more urgency to attend to the population (Figure \ref{e1}-b), the vaccination campaign using $\alpha = 0.98$ can be carried out in fewer days than using $\alpha=0.2$ (86 days / 90 days). 

For $s_2$, the same analysis is performed. Similar results are obtained concerning the prioritization of population groups (see Figure \ref{e1}-c and Figure \ref{e1}-d). Additionally, due to the increase in vaccination capacity in $s_2$, the campaign can be completed in fewer days, finishing on day 66 (see Figure \ref{e1}-d). Therefore, in the two scenarios, the model was able to find a suitable schedule of vaccination in such a way that the expected prioritization was achieved.




To complement the analysis, the vaccination campaign's progress in three out of fourteen macrozones is studied. The proportion of people from groups A and B is greater than 70\% of the corresponding total population. For $s_1$, it can be seen (Figures \ref{e2}-a and \ref{e2}-b) that on day 50 of the vaccination campaign, 80\% of the total population in each of the selected macrozones has been vaccinated. In the same sense, $s_2$ has similar characteristics of attention, where it can be seen (Figures \ref{e2}-c and \ref{e2}-d) that macrozone 9 is attended to with considerable urgency.


\begin{figure}
\centering
\subfloat[]{\includegraphics[width=9cm]{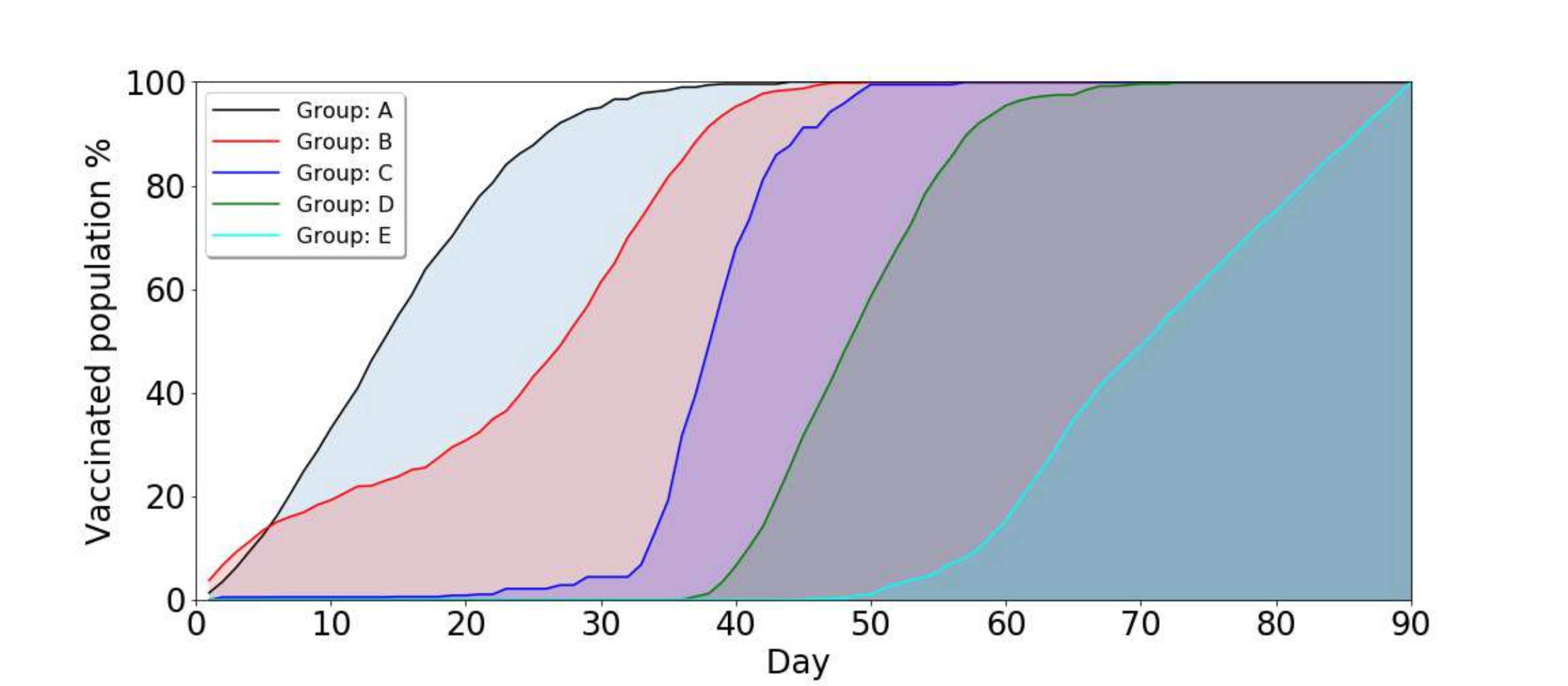}}\hspace{1mm}
\subfloat[]{\includegraphics[width=9cm]{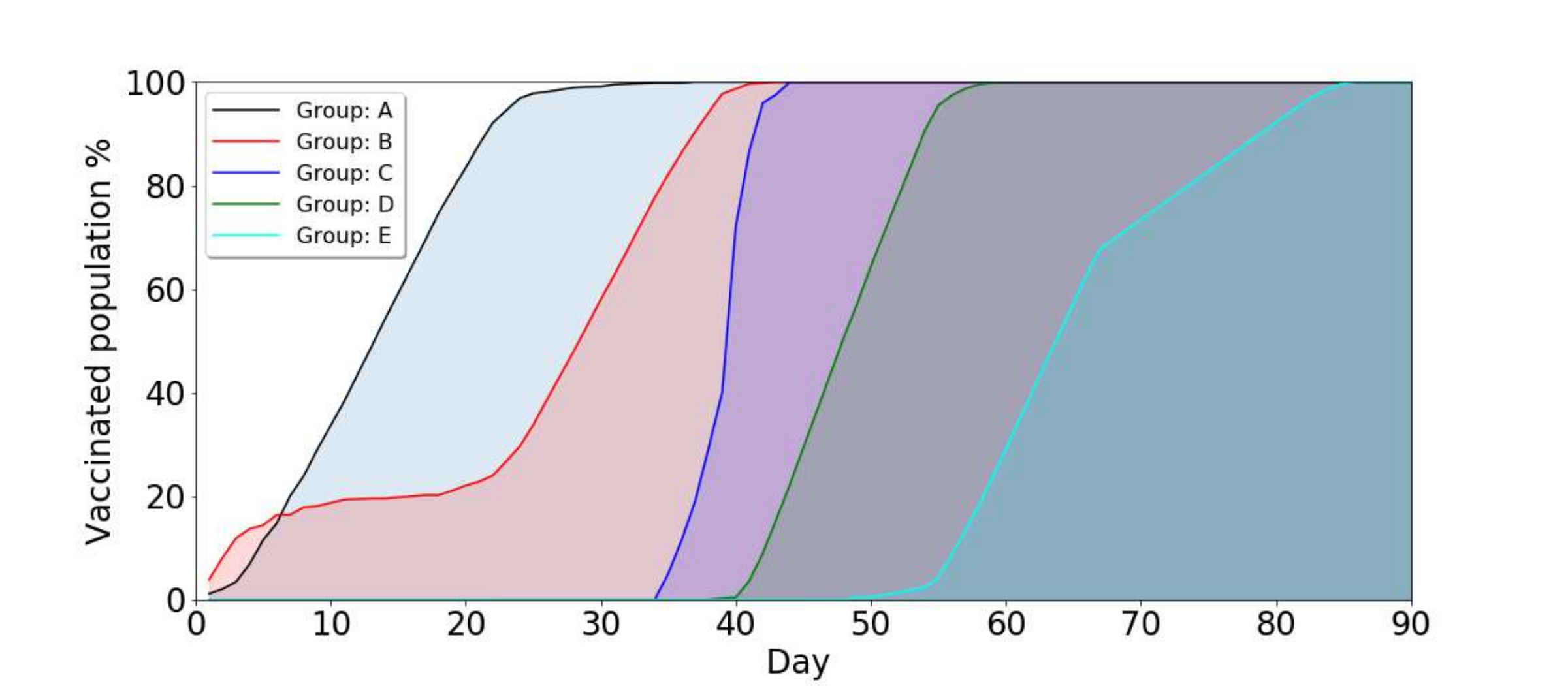}}\vspace{2mm}
\subfloat[]{\includegraphics[width=9cm]{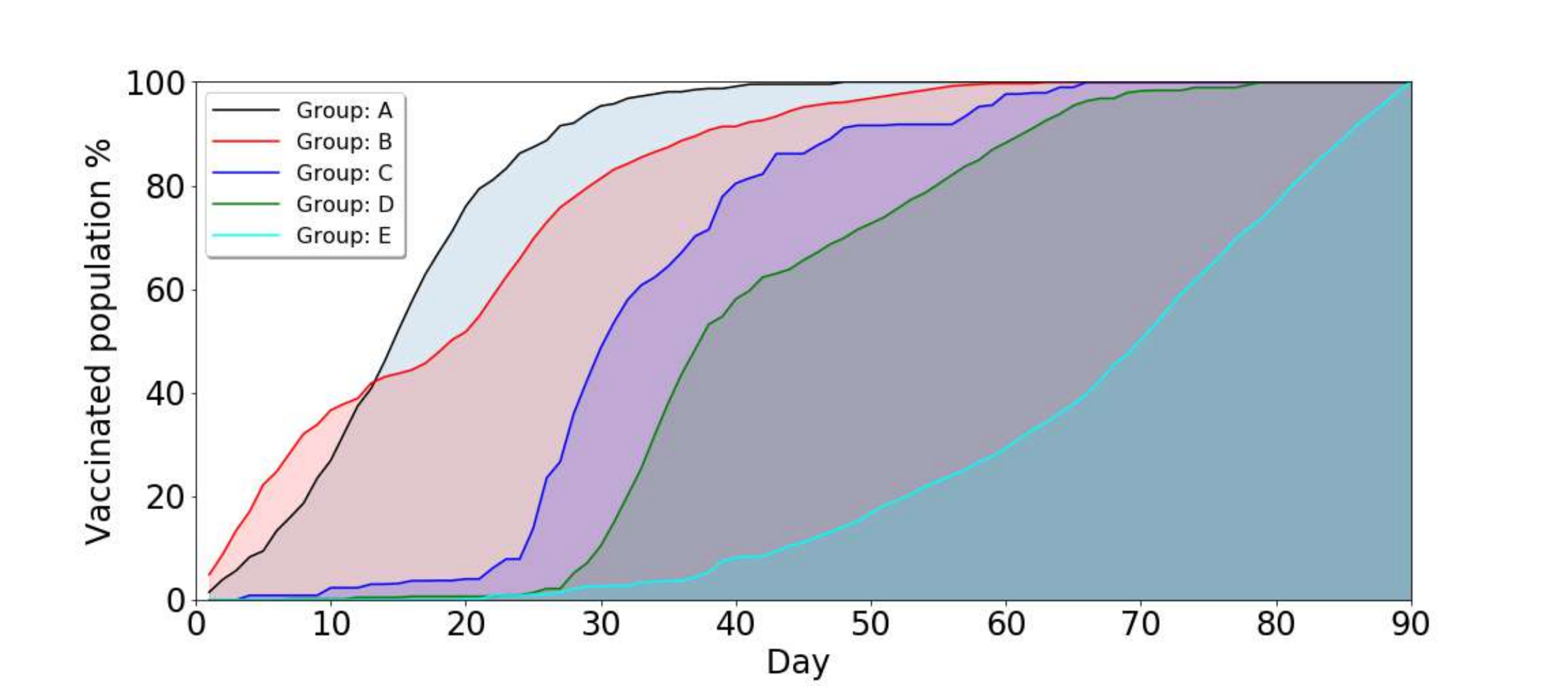}}\hspace{1mm}
\subfloat[]{\includegraphics[width=9cm]{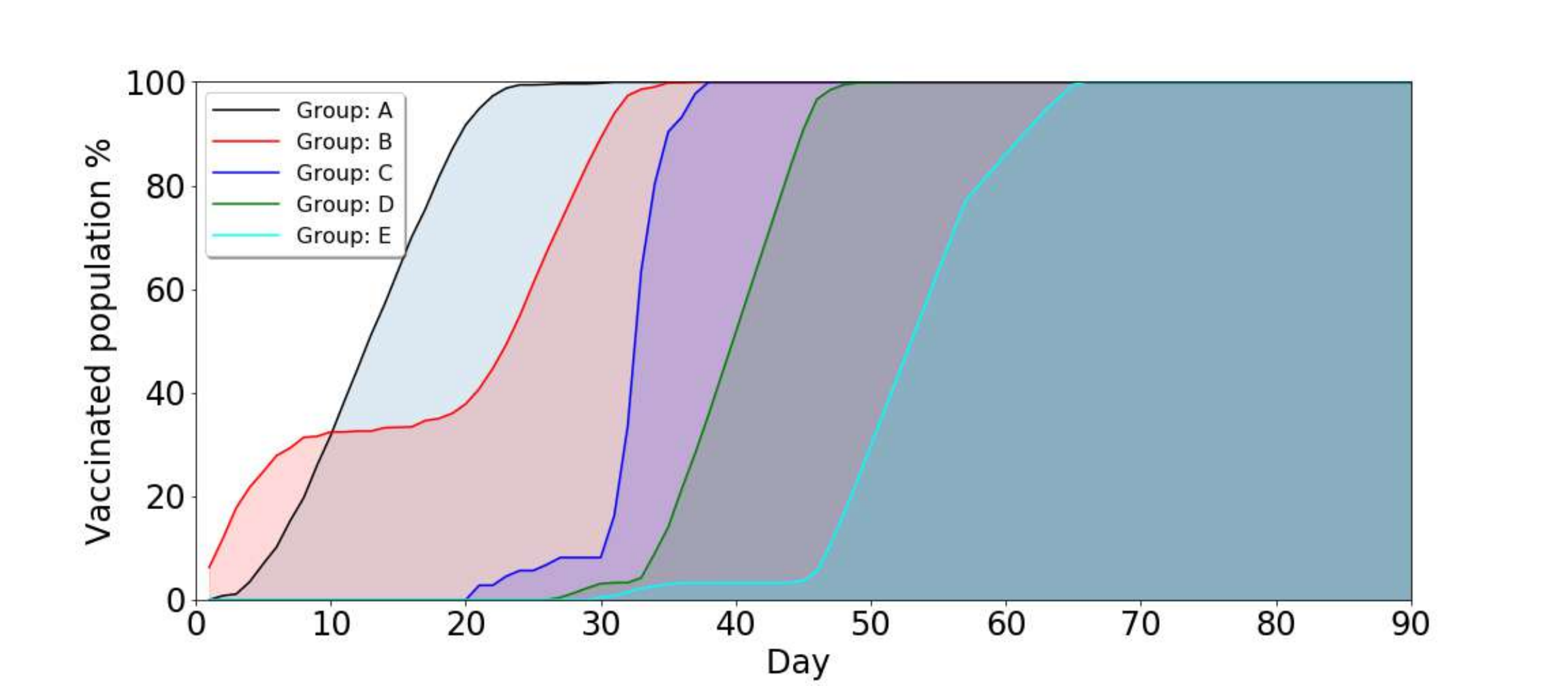}}\vspace{2mm}
\caption{{ \footnotesize Percentage of vaccinated people in each group on day $t$. Scenario $s_1$: (a) when $\alpha=$ 0.2, (b) when $\alpha=$ 0.98. Scenario $S_2$: (c) when $\alpha=$ 0.2, (d) when $\alpha=$ 0.98. }\label{e1}}
\label{e1}
\end{figure}

\begin{figure}
\centering
\subfloat[]{\includegraphics[width=9cm]{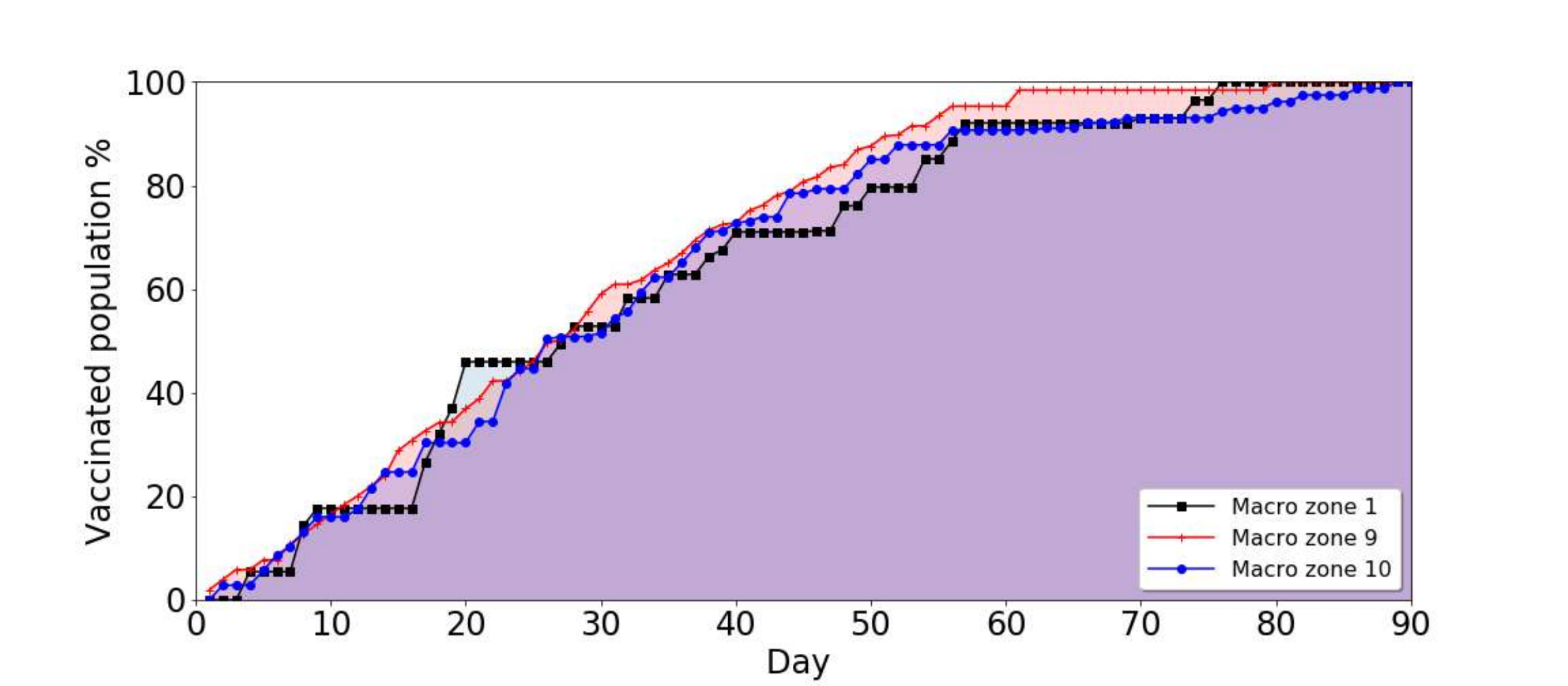}}\hspace{1mm}
\subfloat[]{\includegraphics[width=9cm]{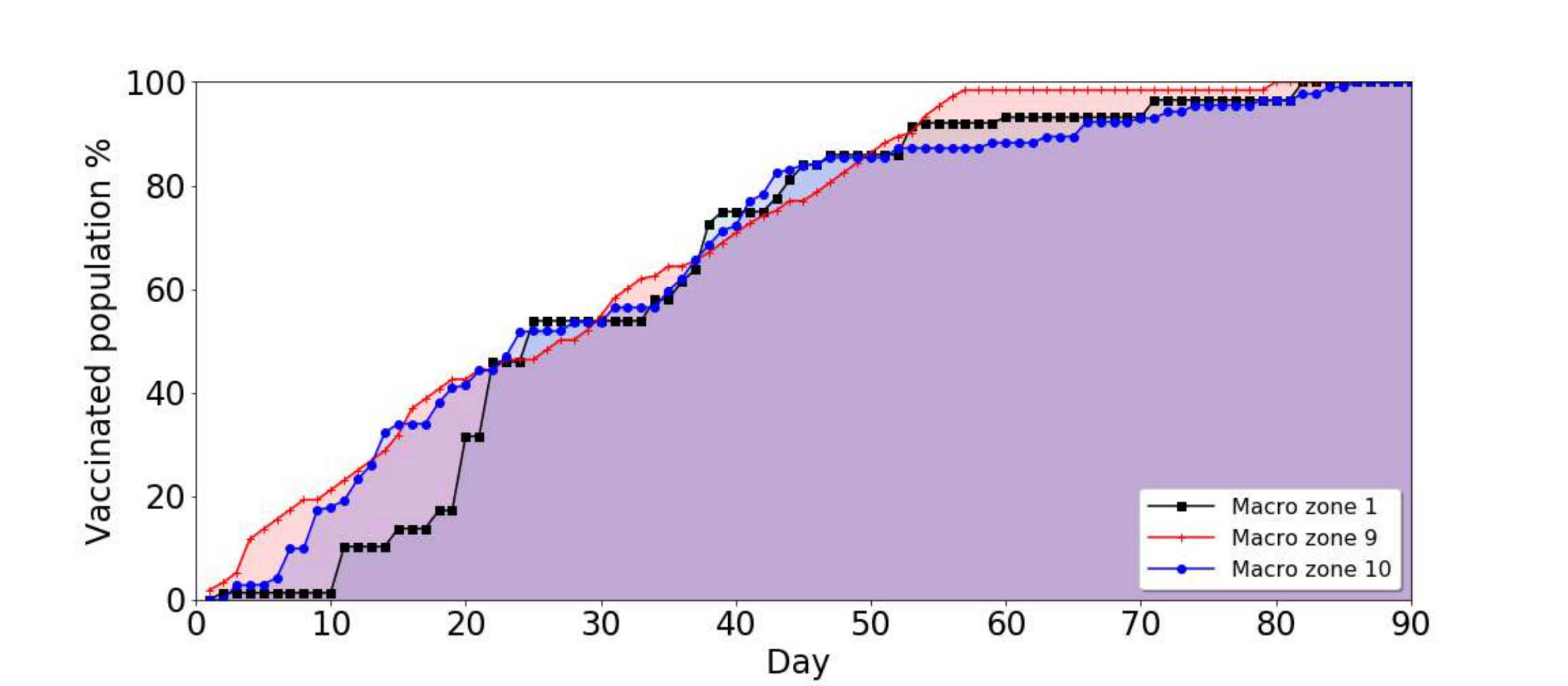}}\vspace{2mm}
\subfloat[]{\includegraphics[width=9cm]{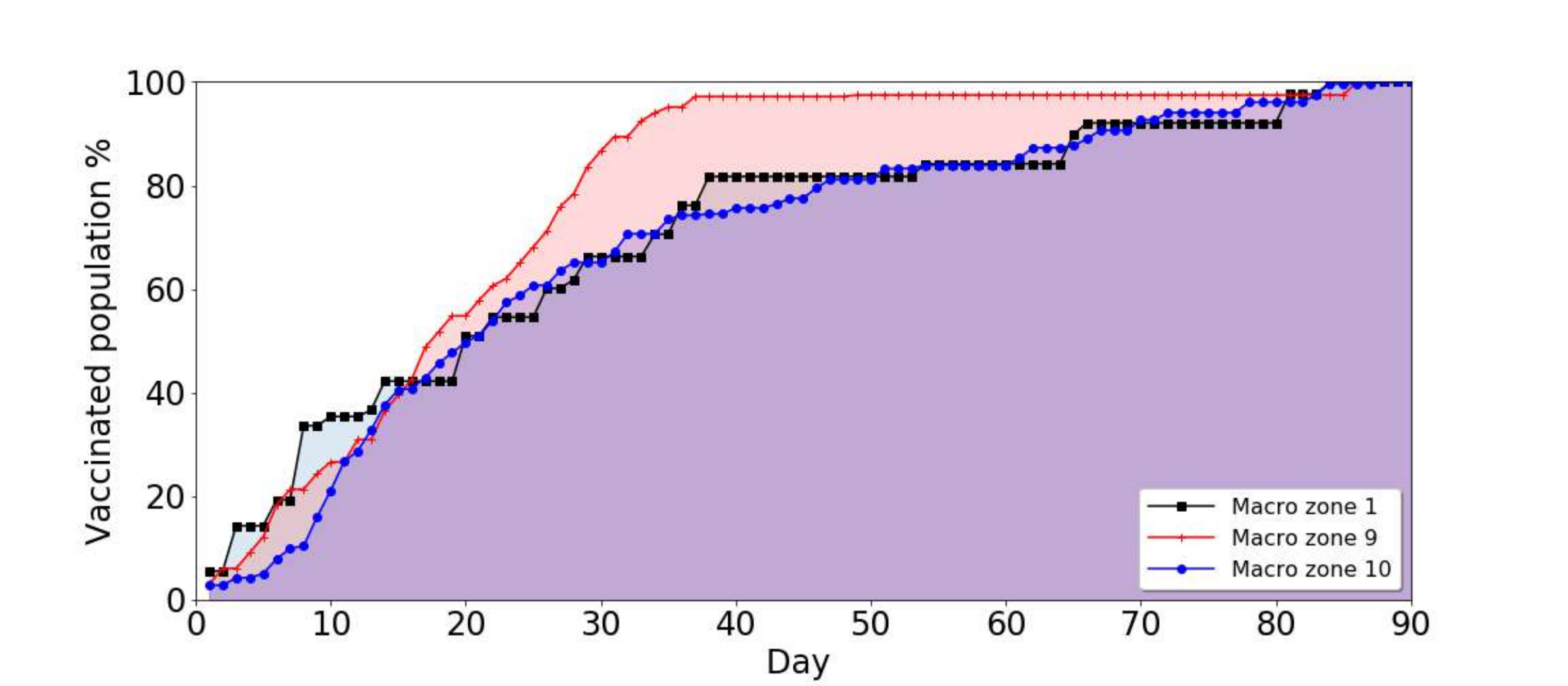}}\hspace{1mm}
\subfloat[]{\includegraphics[width=9cm]{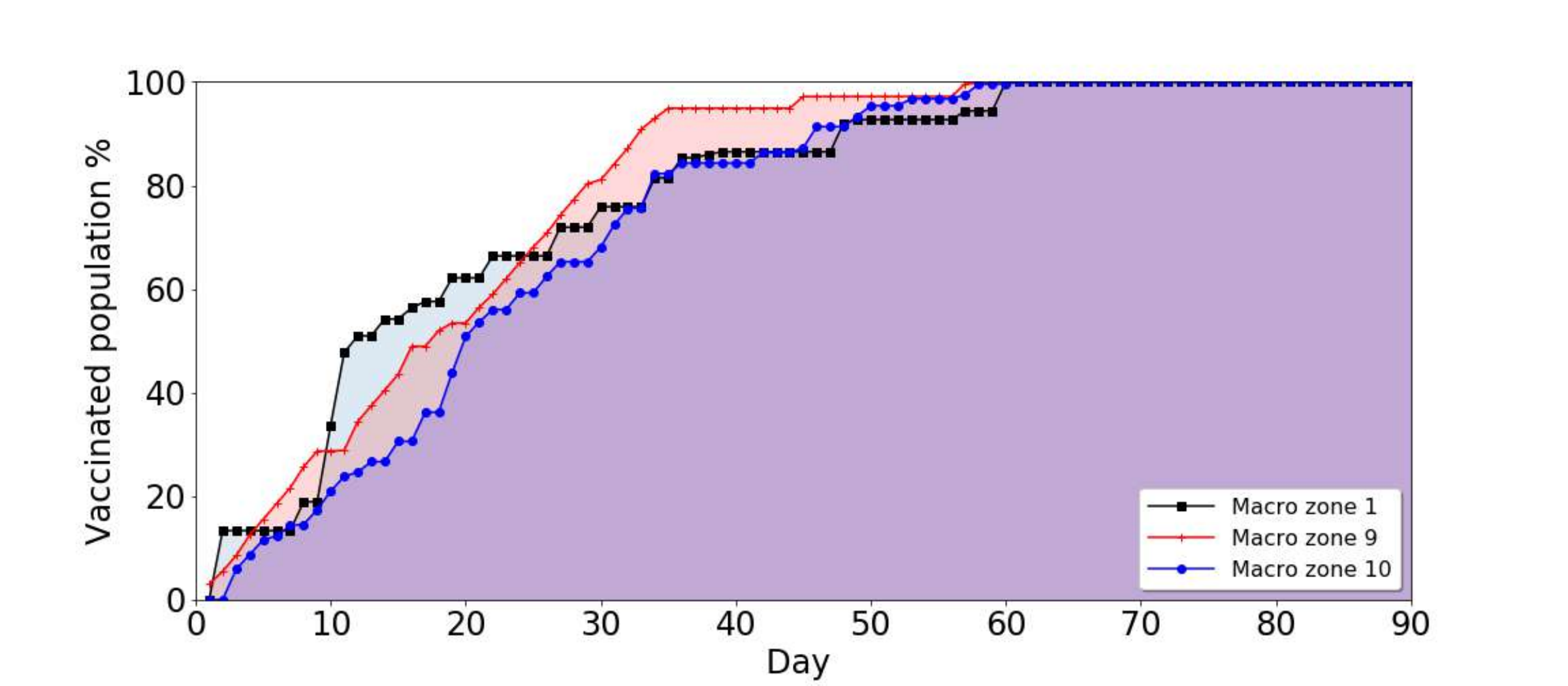}}\vspace{2mm}
\caption{{ \footnotesize Percentage of vaccinated people in each macrozone on day $t$. Scenario $s_1$: (a) when $\alpha=$ 0.2, (b) when $\alpha=$ 0.98. Scenario $s_2$: (c) when $\alpha=$ 0.2, (d) when $\alpha=$ 0.98.}\label{e2}}
\label{e2}
\end{figure}

\subsection{Performance of the weighted bi-objective method}
\label{5.3}


To analyze the bi-objective model's performance respecting the solutions founds, extensive experiments were carried out modifying the $\alpha$ parameter. It is important to note that the scenarios' capacity parameters establish lower and upper limits to both centers' possible use. In particular, given the goal to satisfy the vaccination demand, 115,800 vaccines must be applied. Since permanent centers' total vaccination capacity in the maximum duration of the campaign is $623\times90=56,070$ vaccines, then, with data from $s_1$, at least $59,730$ vaccines must be supported by the temporary centers. Thus, at least the $51.58\%$ of the vaccines must be applied in temporary centers. Similarly, utilizing data from $s_2$, at least $22.28\%$ of the vaccines must be applied in temporary centers. Using the same above analysis to the temporary centers' capacity constraint, upper bounds for the temporary centers' utilization percentage could be calculated. These facts help to a better understanding of the results presented in Table 7.

Table \ref{resumen} presents different indicators for the two scenarios, which help to identify how the vaccination campaign can be carried out, given the resources available. For $s_1$, it can be seen that there is a significant use of temporary centers, 55.2\% for $\alpha= 0.2$. The objective function $f_1^{*}$ decreases slowly when the value of $\alpha$ is incremented while the value of $f_2^{*}$ remain constant, except for the last value $\alpha= 0.98$. This relation between values of $\alpha$ and $f_1^{*}$ are coherent with the definition of $f_1$.
Interestingly, the decrease of $f_1^{*}$ is sufficient to improve the termination of the prioritized groups $A$ and $B$, advancing its completion. For $\alpha= 0.98$, this fact is even clearer, improving not only the first two groups but also finishing the campaign earlier (day 86).
Note that it is possible to reduce the campaign's completion by only a few days due to the available resources.  Since in $s_2$,   there is more capacity in the permanent centers (and also of vaccines) than $s_1$, for $\alpha \geq 0.92$,  campaign ends earlier in $s_2$ than it does in $s_1$, as is detailed below.

On the other hand, an increase in the permanent centers' vaccination capacity allows for a campaign with less dependence on temporary vaccination centers. This can be seen in $s_2$, where it starts at 33.7\% and finishes at 49.2\% of usage. Moreover, the greater use of resources translates into an advance in the vaccination campaign, reaching the end of the campaign 24 days before the maximum day allowed. So, the optimization model is able to provide a broad range of schedules to the decision-maker by managing the $\alpha$ parameter. Finally, we highlight that high running times could be necessary to solve the problem due to the model's complexity. In these experiments, the maximum runtime was 14,400 seconds. Therefore, at least for instances of 70 neighborhoods, ten temporary centers, and a horizon of 90 days, the proposed model can achieve solutions in reasonable times for tactical planning.

\begin{table}
\centering
\footnotesize
\begin{tabular}{|c|c|c|c|c|c|c|}
\hline
Scenario                    & $\alpha$  & $\mathcal{P}$& $f_{1}^{*}$    & $f_{2}^{*}$    & $\mathcal{D}_{A}$-$\mathcal{D}_{B}$-$\mathcal{D}_{C}$-$\mathcal{D}_{D}$-$\mathcal{D}_{E}$ \\ \hline
\multirow{9}{*}{$s_{1}$} & 0.2  & 55.20\% & 0.49756 & 0.15232 & 44-50-57-73-90 \\
                             & 0.4   & 55.30\%  & 0.49418 & 0.15232 & 42-47-51-71-90 \\
                             & 0.6 & 55.60\% & 0.49267 & 0.15232 & 38-45-48-68-90 \\
                             & 0.8 & 55.60\%  & 0.49232 & 0.15232 & 37-46-49-67-90 \\
                             & 0.9 & 55.61\% & 0.49225 & 0.15232 & 33-47-48-66-90 \\
                             & 0.92 & 55.61\% & 0.49210 & 0.15232 & 38-45-49-68-90 \\
                             & 0.94 & 55.61\%  & 0.49201 & 0.15232 & 33-48-48-68-90 \\
                             & 0.96 & 55.61\%  & 0.49197 & 0.15232 & 36-44-48-65-90 \\
                             & 0.98 & 57.69\%  & 0.48994 & 0.23179 & 37-43-44-60-86 \\ \hline
\multirow{9}{*}{$s_{2}$} & 0.2   & 33.70\%  & 0.47752 & 0.20561 & 49-67-79-90-90 \\
                             & 0.4   & 33.70\%  & 0.47692 & 0.20561 & 42-59-73-90-90 \\
                             & 0.6   & 33.70\%  & 0.47690 & 0.20561 & 41-59-61-73-90 \\
                             & 0.8  & 33.70\%  & 0.47686 & 0.20561 & 39-59-61-72-90 \\
                             & 0.9  & 36.00\%&  0.46917 & 0.26168 & 31-47-49-60-90 \\
                             & 0.92 & 39.70\%   & 0.46377 &	0.31464 & 28-44-46-57-84 \\
                             & 0.94  & 43.00\%  & 0.45926 & 0.37383 & 34-41-45-56-78 \\
                             & 0.96  & 45.90\%  & 0.45654 & 0.42679 & 27-38-41-53-72 \\
                             & 0.98 & 49.20\%  & 0.45474 & 0.48598 & 36-37-38-52-66 \\ \hline
\end{tabular}
\caption{\footnotesize Results for different $\alpha$ values. $\mathcal{P}$ represents the percentage of patients vaccinated at temporary centers. 
$f_{1}^{*}$ and $f_{2}^{*}$ are the values of the objective functions already normalized. $\mathcal{D}_{p}$ represents the day on which the risk group $p$ $\in$ $\{A,B,C,D,E\}$ has finished being vaccinated. }
\label{resumen}
\end{table}

\subsection{Temporary centers covering}
\label{5.4}
This section presents the deployment of the vaccination campaign in different macrozones, through the care provided by permanent and temporary centers. To analyze how each of the macrozones is cared for, the distribution of the population presented in Figure \ref{fy} is considered, and the results of scenario $s_2$ with $\alpha=0.98$ are illustrated in Figure \ref{em2}.

One macrozone with great demand is No. 6. For this reason, it must be attended to with great assistance from permanent and temporary centers. In other sectors, the model prioritizes the use of one of the vaccination centers over another. In the case of the sectors which have to use permanent centers to a greater extent, there are macrozones 9, 4, and 5. The deployment of temporary centers in sectors 8 and 14 stand out because, after macrozone 6, these are the most populous sectors, which account for 22\% of the population in the San Bernardo region. Note that macrozones 8 and 14 include many people from groups $A$ and $B$. Also, note that macrozone 12 is the smallest macrozone and its intensity in the figure is minimum. The results for scenario $s_1$ are analogous to those of scenario $s_2$.



\begin{figure}
\centering
\subfloat[]{\includegraphics[width=4cm]{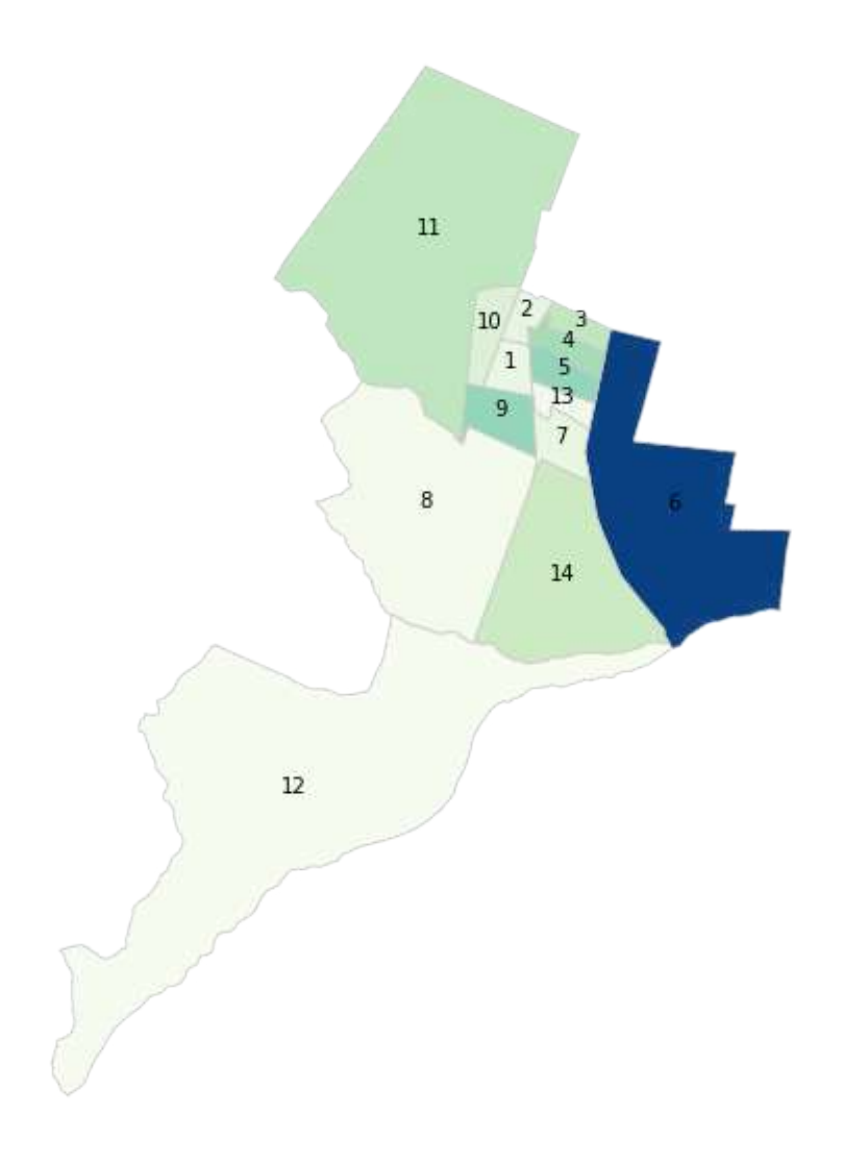}}\hspace{1mm}
\subfloat[]{\includegraphics[width=5.4cm]{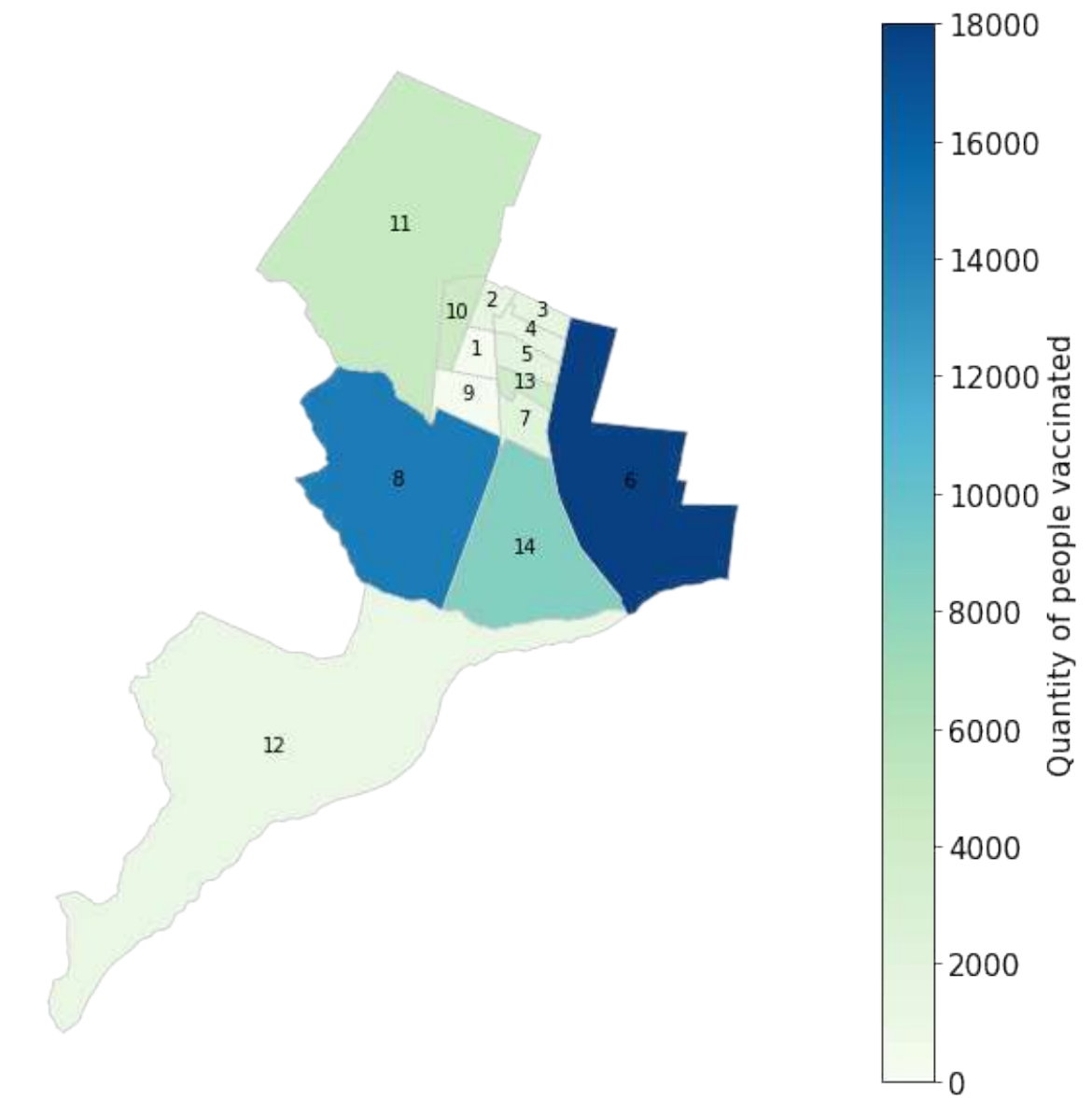}}\vspace{2mm}
\caption{{\footnotesize Scenario $s_2$. Deployment of the vaccination campaign in the macrozones: (a) permanent centers and (b) temporary centers, both at $\alpha =0.98$.}\label{em2}}
\label{em2}
\end{figure}


\subsection{Level of risk and the use of temporary centers}
\label{5.5}

The population groups $A$ and $B$ defined in this vaccination campaign have mobility difficulties since they include people with advanced age and chronic diseases. This means that these groups' travel to permanent centers can be complicated or even impossible. In this way, the use of temporary centers becomes relevant because it can make vaccination more accessible to these people. Table \ref{tab:temp2} shows two metrics, which allow identifying the use of temporary centers by population groups. Firstly, the metric $\mathcal{T}_p$ is calculated as presented in equation (\ref{eq:rt}), it represents the percentage of each group ($p \in P$) of the total that was vaccinated with temporary centers. Secondly, the metric $\mathcal{G}_p$ presented in equation (\ref{eq:rg}) allows us to analyze what percentage of the total ($pop_{pl}$) of each group was vaccinated in temporary centers. It is important to highlight that regarding both scenarios and values of $\alpha$, people from groups $A$ and $B$ achieve 59\% and 71\%, respectively, from the total vaccinated in temporary centers (metric $\mathcal{T}_p$).  

\begin{align}
     \mathcal{T}_p&= \sum_{l \in \left\lbrace L_{k}: k \in K \right\rbrace, t \in T, j \in J}\dfrac{\gamma_{lptj}}{\sum_{p \in P}\gamma_{lptj}}\times 100 && \forall p \in P \label{eq:rt}\\
    \mathcal{G}_p&= \sum_{l \in \left\lbrace L_{k}: k \in K \right\rbrace, t \in T, j \in J}\dfrac{\gamma_{lptj}}{\sum_{p \in P}pop_{pl}}\times 100 && \forall p \in P \label{eq:rg}
\end{align}

Furthermore, when we look at the case of scenario $s_1$ through metric $\mathcal{G}_p$, it can be seen that, for each group, when taking a plan with $\alpha=0.2$,  between 41\% and 63\%  of the people are served by this type of resource. This is since to carry out the campaign, the capacity of temporary centers is greater than the permanent centers. On the other hand, in the case of scenario, $s_2$ and also using the metric $\mathcal{G}_p$, because of the greater capacity of permanent centers, the use of temporary centers was reduced but still achieves, in average,
$56\%$ for $\alpha=0.2$ and $47\%$ for $\alpha=0.98$ 
of the total. Finally, in the scenario $s_1$, note that in both metrics, the results are similar, regardless of the value of $\alpha$. Contrarily, in scenario $s_2$,  the results are different depending on the value of alpha for metric $\mathcal{G}_p$. Specifically, for all the groups, vaccination using metric $\mathcal{G}_p$ include clearly more people than those vaccinated using metric $\mathcal{T}_p$. This difference is clearer in groups A, B, C, and D than in group E. 
 This effect comes from the high value of $\alpha$ implies that the second objective has a lower weight. Consequently, it is possible to spend resources on the temporary centers independently of the costs and only limited by the vaccines available.
Therefore, the model can assign people to both types of centers according to the centers' capacity and weights defined in the objective functions.



\begin{table}
\centering
\footnotesize
\begin{tabular}{ccc|ccccc|c}
\cline{4-8}
                               &                                            &         & \multicolumn{5}{c|}{Population group} &                           \\ \hline
\multicolumn{1}{|c|}{Scenario} & \multicolumn{1}{c|}{$\alpha$}                 & Metrics & A     & B     & C     & D     & E     & \multicolumn{1}{c|}{Mean} \\ \hline
\multicolumn{1}{|c|}{\multirow{4}{*}{$s_1$}} & \multicolumn{1}{c|}{\multirow{2}{*}{0.2}} & $\mathcal{T}_p$ & 28\% & 31\% & 3\% & 21\% & 17\% & \multicolumn{1}{c|}{20\%} \\
\multicolumn{1}{|c|}{}         & \multicolumn{1}{c|}{}                      & $\mathcal{G}_p$      & 58\%  & 63\%  & 50\%  & 58\%  & 41\%  & \multicolumn{1}{c|}{54\%} \\ \cline{2-9} 
\multicolumn{1}{|c|}{}         & \multicolumn{1}{c|}{\multirow{2}{*}{0.98}} & $\mathcal{T}_p$       & 27\%  & 31\%  & 3\%   & 22\%  & 17\%  & \multicolumn{1}{c|}{20\%} \\
\multicolumn{1}{|c|}{}         & \multicolumn{1}{c|}{}                      & $\mathcal{G}_p$      & 59\%  & 64\%  & 48\%  & 64\%  & 44\%  & \multicolumn{1}{c|}{56\%} \\ \hline
\multicolumn{1}{|c|}{\multirow{4}{*}{$s_2$}} & \multicolumn{1}{c|}{\multirow{2}{*}{0.2}} & $\mathcal{T}_p$ & 31\% & 40\% & 3\% & 23\% & 3\%  & \multicolumn{1}{c|}{20\%} \\
\multicolumn{1}{|c|}{}         & \multicolumn{1}{c|}{}                      & $\mathcal{G}_p$       & 59\%  & 64\%  & 48\%  & 64\%  & 44\%  & \multicolumn{1}{c|}{56\%} \\ \cline{2-9} 
\multicolumn{1}{|c|}{}         & \multicolumn{1}{c|}{\multirow{2}{*}{0.98}} & $\mathcal{T}_p$      & 27\%  & 29\%  & 2\%   & 23\%  & 19\%  & \multicolumn{1}{c|}{20\%} \\
\multicolumn{1}{|c|}{}         & \multicolumn{1}{c|}{}                      & $\mathcal{G}_p$       & 50\%  & 52\%  & 34\%  & 56\%  & 42\%  & \multicolumn{1}{c|}{47\%} \\ \hline
\end{tabular}
\caption{\footnotesize Distribution of the groups that were vaccinated in temporary centers for each scenario.}
\label{tab:temp2}
\end{table}

\begin{figure}
\centering
\subfloat[]{\includegraphics[width=11.5cm]{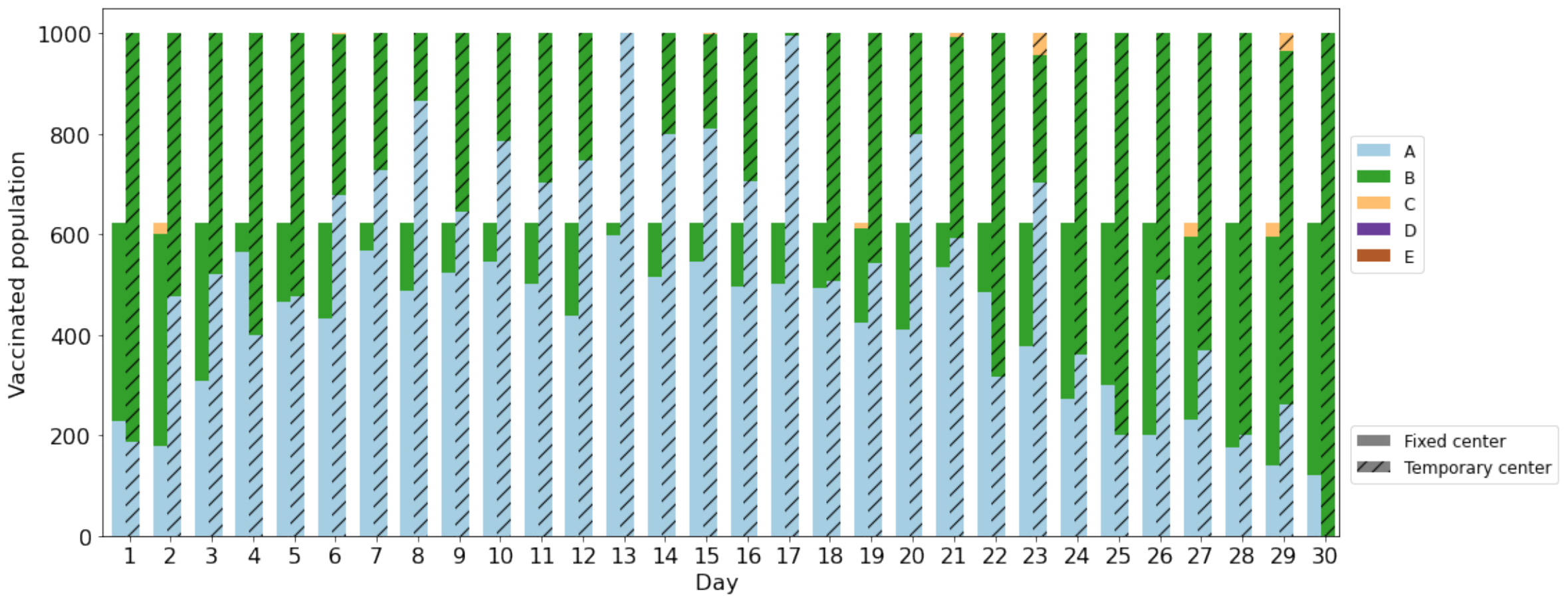}}\hspace{1mm}
\subfloat[]{\includegraphics[width=11.5cm]{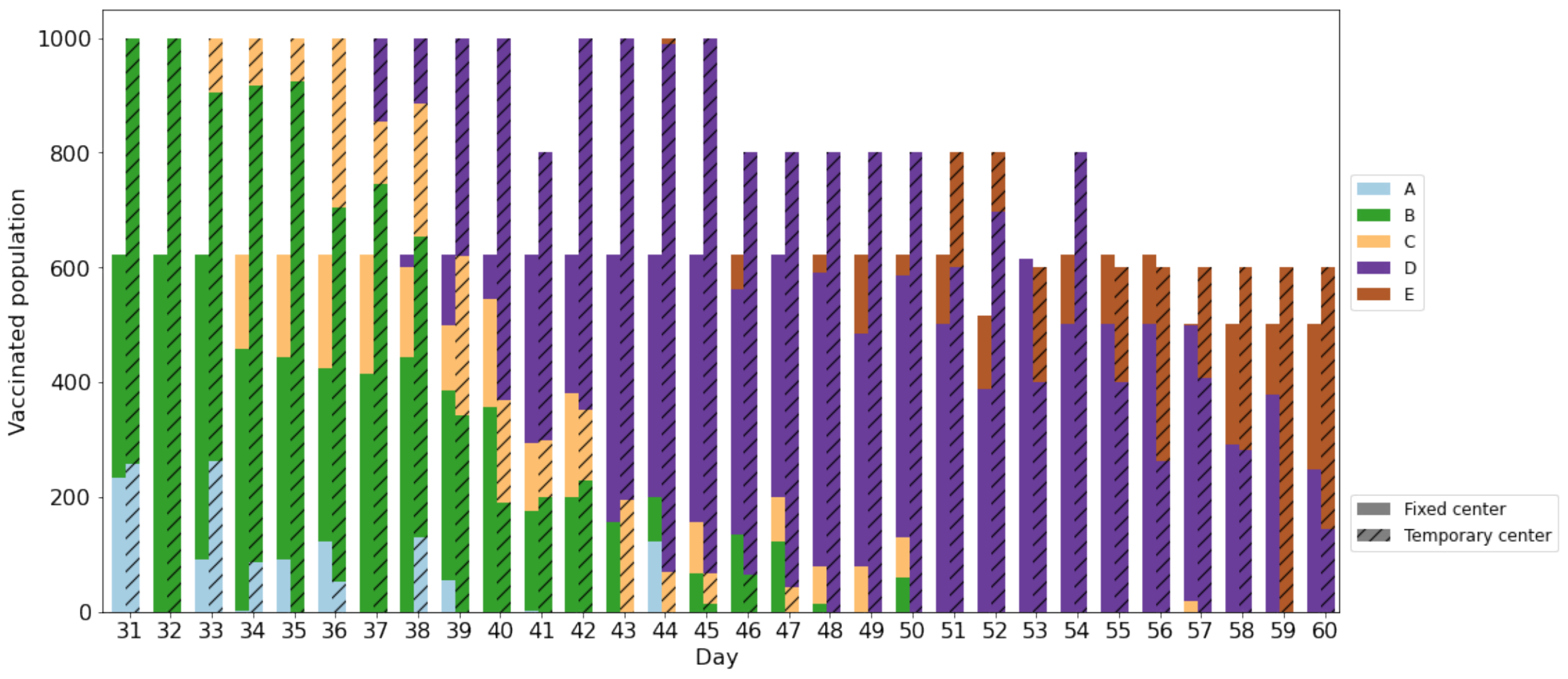}}\vspace{2mm}
\subfloat[]{\includegraphics[width=11.5cm]{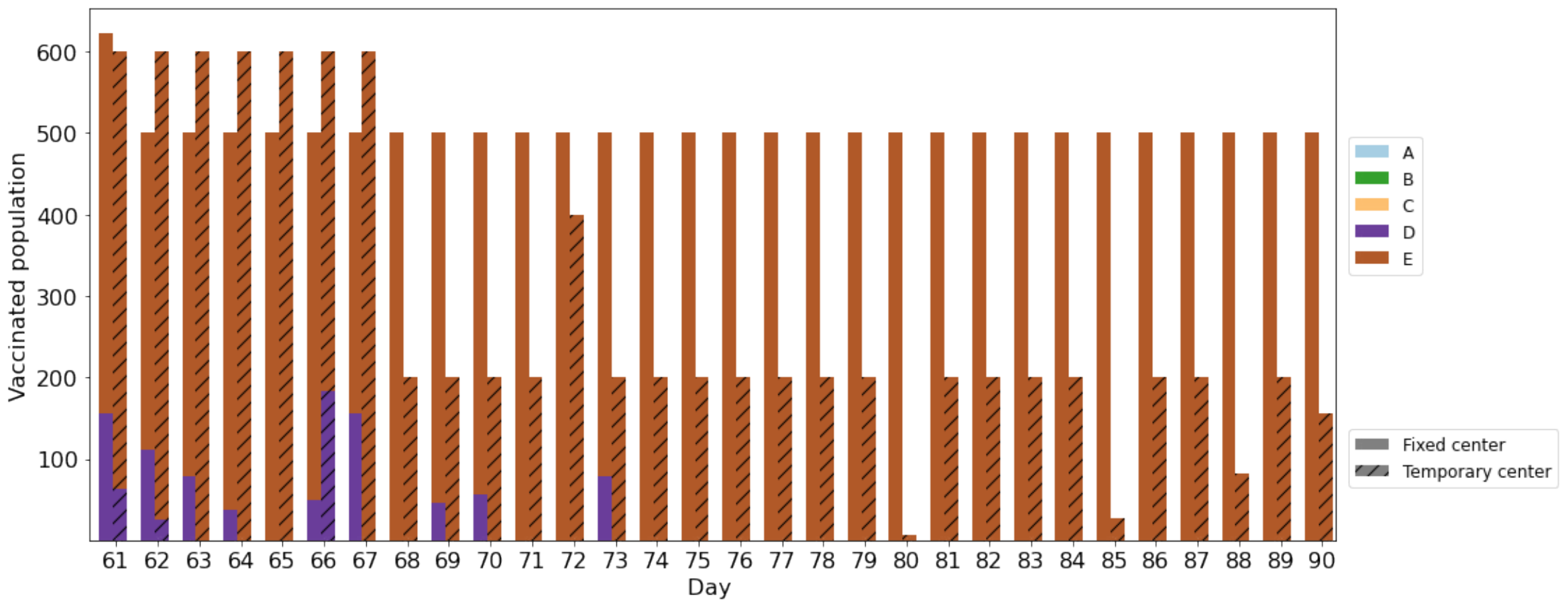}}\vspace{2mm}
\caption{{ \footnotesize Progress of the vaccination campaign differentiated by population group for permanent and temporary centers. Scenario $s_1$ with $\alpha=0.2$.  }\label{barra_e1_0.2}}
\end{figure}
\begin{figure}
\centering
\subfloat[]{\includegraphics[width=11.5cm]{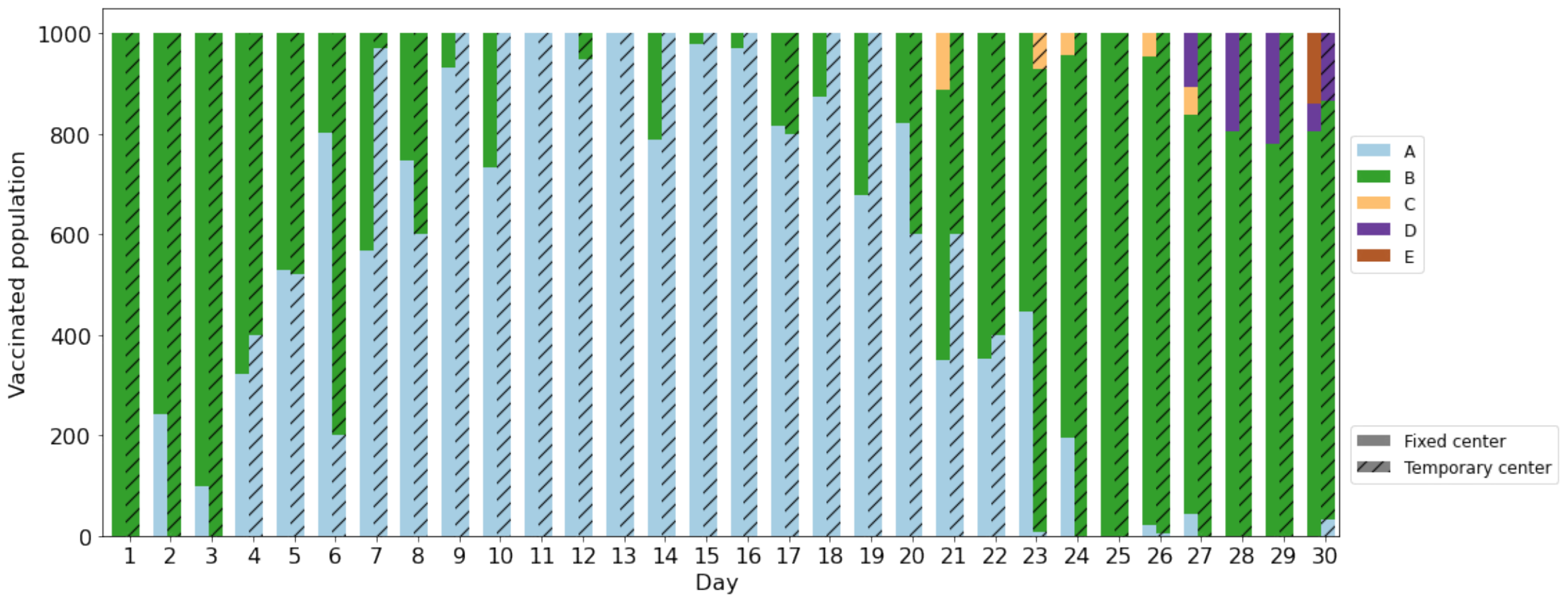}}\hspace{1mm}
\subfloat[]{\includegraphics[width=11.5cm]{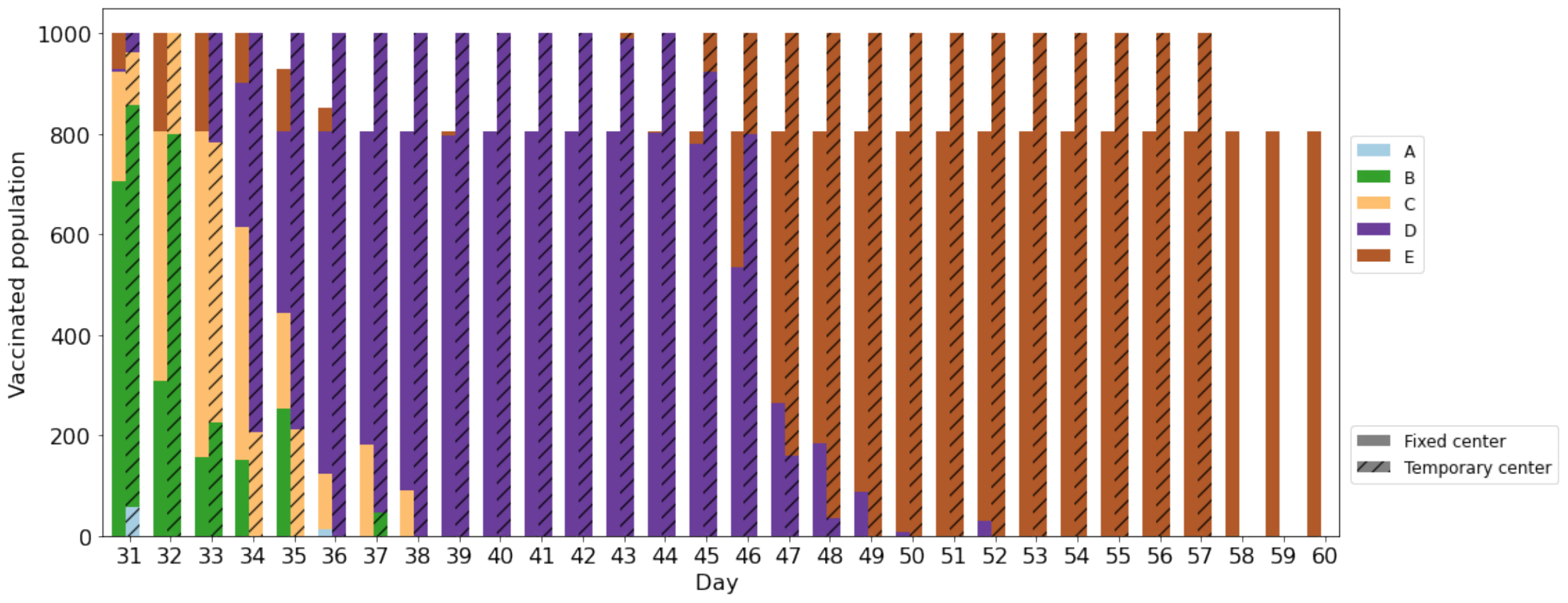}}\vspace{2mm}
\subfloat[]{\includegraphics[width=11.5cm]{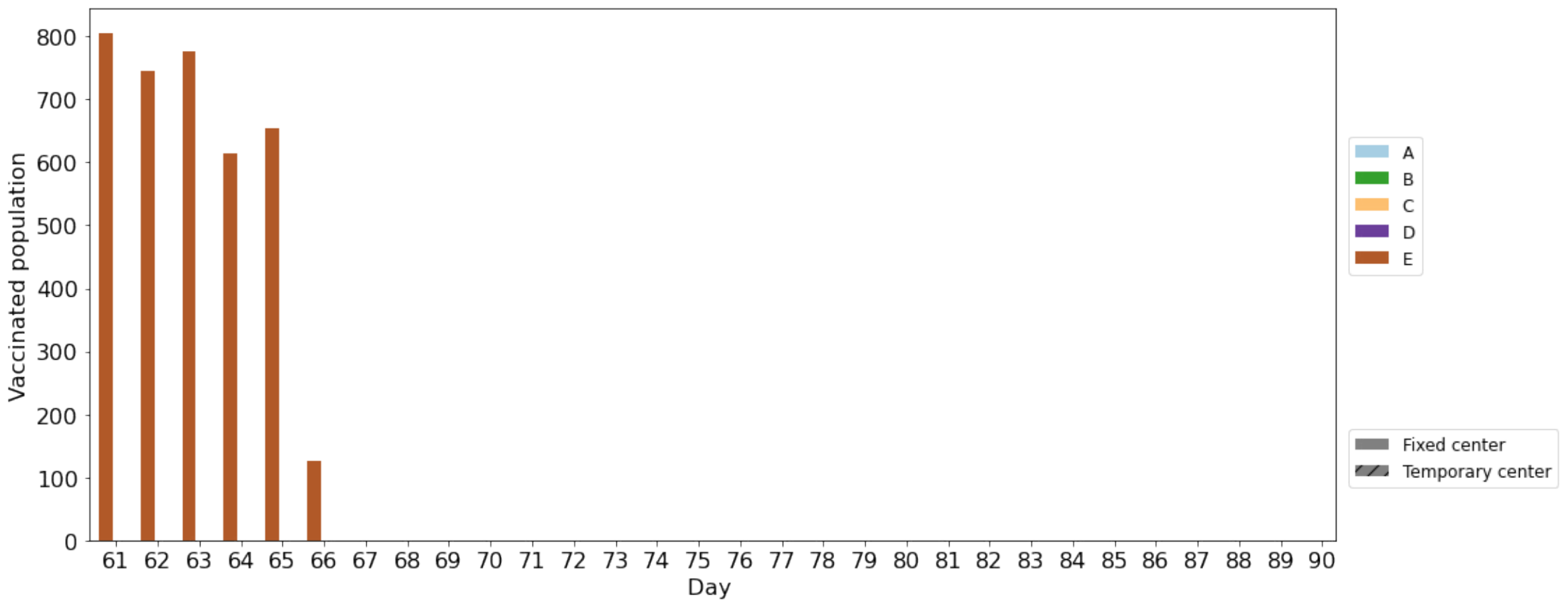}}\vspace{2mm}
\caption{{ \footnotesize Progress of the vaccination campaign differentiated by population group for permanent and temporary centers. Scenario $s_2$ with $\alpha=0.98$.}\label{barra_e2_0.98}}
\end{figure}

Another relevant aspect is to visualize how these population groups' demand is distributed among permanent and temporary centers over time. Figures \ref{barra_e1_0.2} and \ref{barra_e2_0.98} are proposed for this purpose, choosing one value of $\alpha$ for each of the scenarios. The bars represent the number of people vaccinated, the color indicates the population group served and the stripes differentiate between permanent and temporary centers. For the first scenario $s_1$, Figure \ref{barra_e1_0.2} presents an operational plan for $\alpha= 0.2$. It can be seen that temporary centers are intensively used until day 45 (see Fig. \ref{barra_e1_0.2} (a) and (b)). After day 45, their use decreases progressively, demonstrating their relevance to carrying out the campaign. Following the analyses carried out in the previous sections, the population groups' prioritization can be seen. 

Concerning scenario $s_2$, the operating plan was analyzed with $\alpha=0.98$. The results are presented in Figure \ref{barra_e2_0.98}, which shows the increase in vaccination capacity at permanent centers, equaling temporary centers' maximum capacity. Additionally, it is highlighted that temporary vaccination centers use their full capacity until day 57, allowing the duration of the vaccination campaign to be reduced by 23 days. Figures \ref{barra_e1_0.2} and \ref{barra_e2_0.98} illustrate the vaccination scheduling of the groups according to both scenarios.




\subsubsection{Detailed analysis of the use of temporary centers }
\label{5.5.1}

To give a more detailed view of how temporary centers operate, two examples of five-day operating plans for scenario $s_1$ with $\alpha=0.2$ are presented in Table \ref{tabla_m2}. For the vaccination plan of temporary center $j=1$, presented on the top of the Table, the first five days of vaccination were chosen. It can be seen that it serves people who belong to the populations' group $P=\{1,2\}$. On the other hand, looking at the deployment of the temporary center $j=2$ for days 40-44 of the vaccination campaign, it can be seen that it serves population groups $P=\{3,4\}$. This shows the prioritization of the population groups as the vaccination campaign progresses. In addition, the flexibility of the model to plan vaccination by prioritizing population groups according to their respective risk, as well as the location activities of the temporary centers, allowing vaccination both in the neighborhood in which they are installed and in those that meet the coverage criteria, is also noteworthy. Finally, it is recognized that both vaccination centers were working at maximum capacity during the days analyzed.

\begin{table}
\footnotesize
\centering
\begin{tabular}{c|ccccccc|c|c|c}
\cline{2-10}
                          & \multicolumn{7}{c|}{Neighborhoods, $j=1$} & \multicolumn{2}{c|}{Groups} &                            \\ \hline
\multicolumn{1}{|c|}{Day} & 12       & 13       & 14      & 18       & 26      & 45       & 48       & 1             & 2           & \multicolumn{1}{c|}{Total} \\ \hline
\multicolumn{1}{|c|}{1}   &          &          &         &          &         &          & 200      & 200           &             & \multicolumn{1}{c|}{200}   \\
\multicolumn{1}{|c|}{2}   &          & 200      &         &          &         &          &          & 200           &             & \multicolumn{1}{c|}{200}   \\
\multicolumn{1}{|c|}{3}   &          &          &         & 159      & 41      &          &          & 200           &             & \multicolumn{1}{c|}{200}   \\
\multicolumn{1}{|c|}{4}   &          &          &         &          &         & 200      &          & 200           &             & \multicolumn{1}{c|}{200}   \\
\multicolumn{1}{|c|}{5}   & 131      &          & 69      &          &         &          &          & 131           & 69          & \multicolumn{1}{c|}{200}   \\ \hline

                          & \multicolumn{7}{c|}{Neighborhoods $j=2$} & \multicolumn{2}{c|}{Groups} &                            \\ \hline
\multicolumn{1}{|c|}{Day} & 36   & 37   & 40  & 60  & 61 & 68  & & 3            & 4            & \multicolumn{1}{c|}{Total} \\ \hline
\multicolumn{1}{|c|}{40}   &       & 131  & 69  &     &    &    &  & 69           & 131          & \multicolumn{1}{c|}{200}   \\
\multicolumn{1}{|c|}{41}   & 200  &      &     &     &    &     &    &           & 200          & \multicolumn{1}{c|}{200}   \\
\multicolumn{1}{|c|}{42}   &      &      &     & 110 & 90 &   &   & 42           & 158          & \multicolumn{1}{c|}{200}   \\
\multicolumn{1}{|c|}{43}   &      &      & 200 &     &    &     &         &      & 200          & \multicolumn{1}{c|}{200}   \\
\multicolumn{1}{|c|}{44}   &      &      &     &     &    & 200 &           &    & 200          & \multicolumn{1}{c|}{200}   \\ \hline
\end{tabular}
\caption{Number of vaccines delivered daily for different neighborhoods through the temporary center, $j=1$ and $j=2$.}
\label{tabla_m2}
\end{table}

\subsection{Minimizing the chances of contagion for higher priority groups}
\label{5.6}

To avoid crowding, the temporary centers are beneficial, helping to relieve congestion in permanent facilities and serving those with mobility difficulties and with greater risk of complications upon infection. To deal with the problem presented by these risk groups, constraint (\ref{r102}) can be added, defining a set $U$ for those groups who do not want to be vaccinated in permanent centers.  
To illustrate the impact of this modification, scenario $s_2$ was considered, and we chose group A as the one which must be attended to only in temporary centers. The result of the model including constraint (\ref{r102}) can be seen in Figure \ref{figcip}. Here we can note that is possible to bring the end date of the vaccination campaign even farther forward in this case (for $\alpha=0.98$) while the prioritization of the groups is still preserved. (Figure \ref{figcip}-b).

\begin{align}
    \sum_{l \in L_{k}}\phi_{lpti} &= 0 && \forall p \in U, t \in T, i \in I  \label{r102}
\end{align}

\begin{figure}
\centering
\subfloat[]{\includegraphics[width=7.5cm]{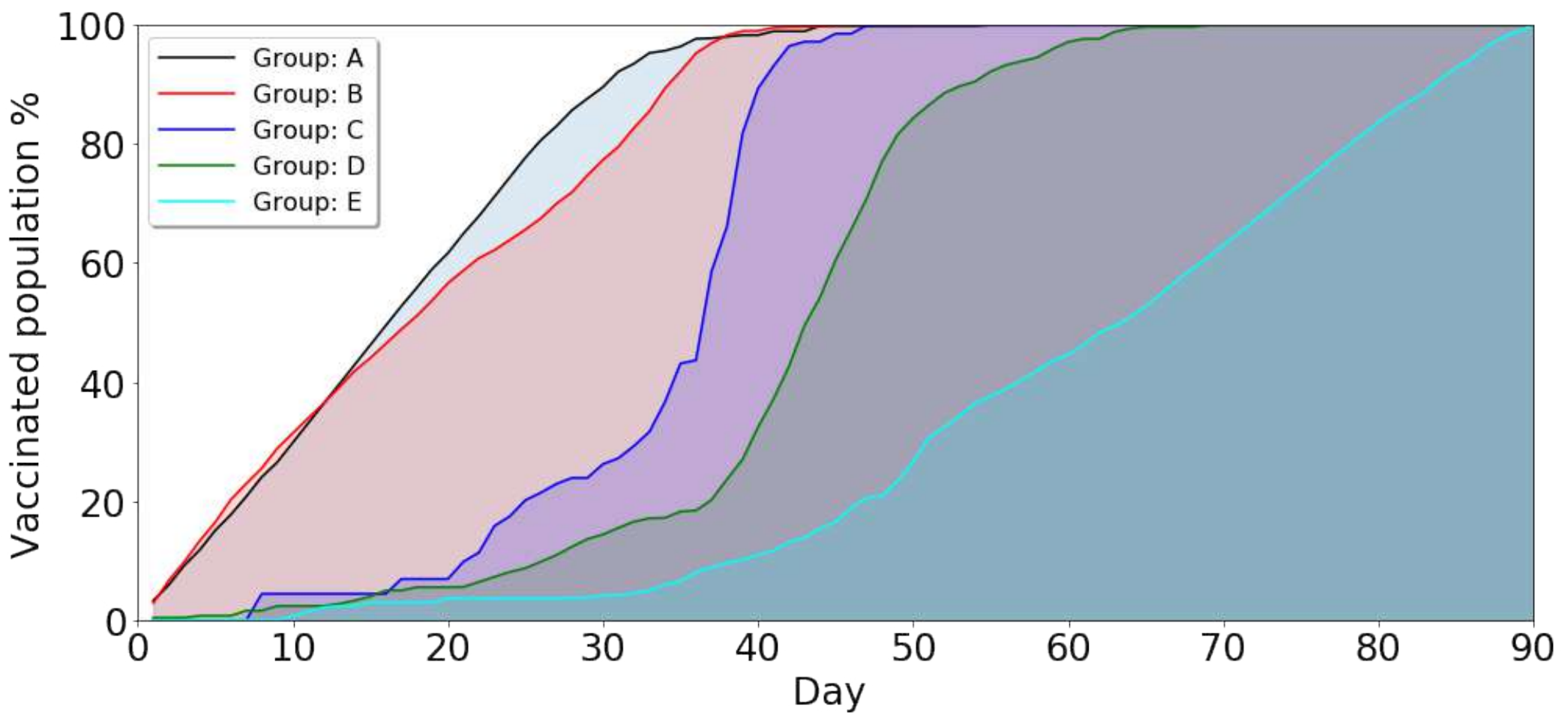}}\hspace{0mm}
\subfloat []{\includegraphics[width=7.5cm]{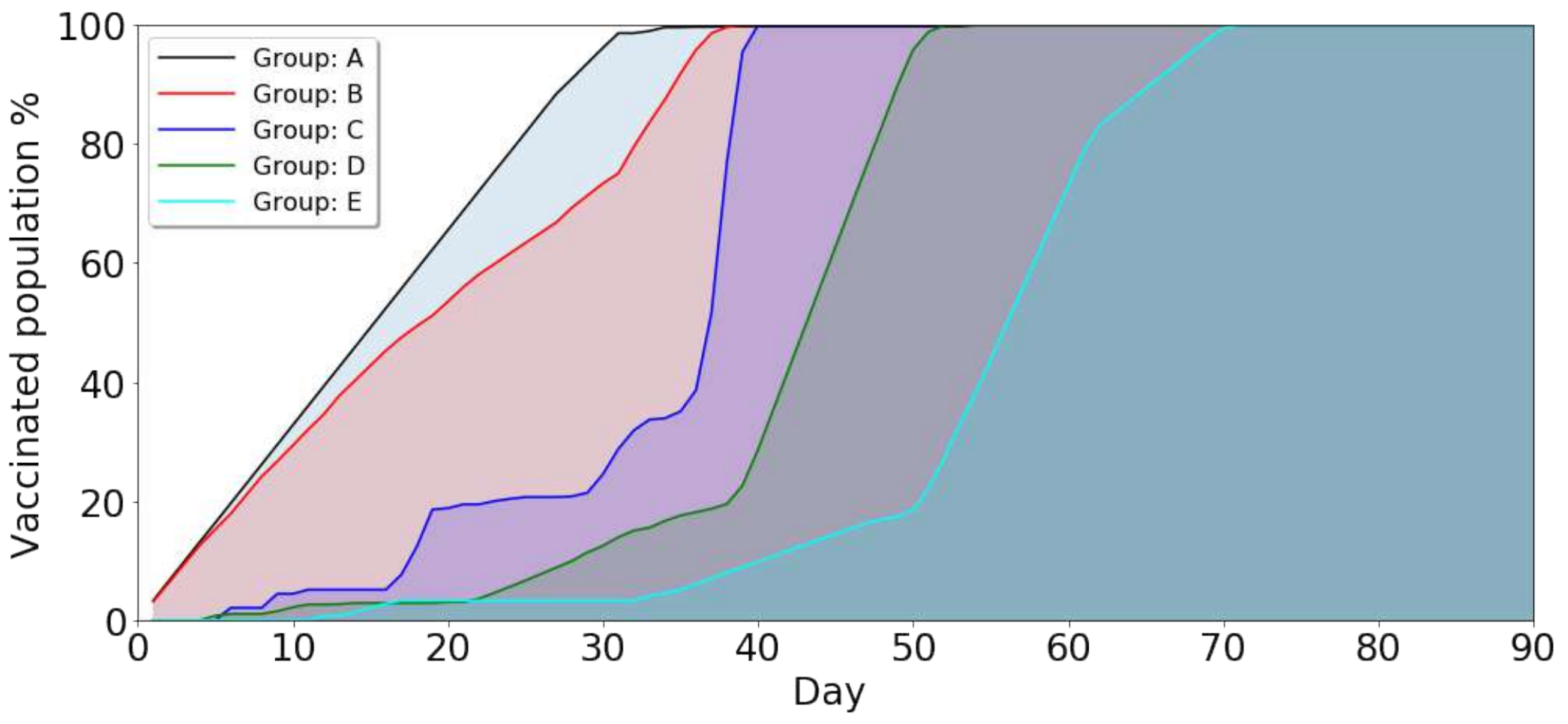}}\vspace{1mm}
\caption{{ \footnotesize Percentage of vaccinated people of each group at day $t$ (a) $\alpha=$ 0.2, (b) $\alpha=$ 0.98. Scenario $s_2$.  
}\label{figcip}}
\label{figcip}
\end{figure}


\begin{table}
\footnotesize
\centering
\begin{tabular}{cc|ccccc|}
\cline{3-7}
                                            &         & \multicolumn{5}{c|}{Population Group} \\ \hline
\multicolumn{1}{|c|}{$\alpha$}                 & Metrics & A      & B     & C     & D     & E    \\ \hline
\multicolumn{1}{|c|}{\multirow{2}{*}{0.2}}  & $\mathcal{T}_p$       & 59\%   & 12\%  & 3\%   & 20\%  & 6\%  \\
\multicolumn{1}{|c|}{}                      & $\mathcal{G}_p$       & 100\%  & 20\%  & 35\%  & 45\%  & 12\% \\ \hline
\multicolumn{1}{|c|}{\multirow{2}{*}{0.98}} & $\mathcal{T}_p$       & 50\%   & 11\%  & 2\%   & 20\%  & 18\% \\
\multicolumn{1}{|c|}{}                      & $\mathcal{G}_p$       & 100\%  & 21\%  & 34\%  & 53\%  & 41\% \\ \hline
\end{tabular}
\caption{\footnotesize Distribution of the groups that were vaccinated in temporary centers. Scenario $s_2$.}
\label{mob_cip}
\end{table}

The distribution of the groups served with temporary centers, measured by the metrics $\mathcal{T}_p$ and $\mathcal{G}_p$ (presented in the section \ref{5.5}) is shown in Table \ref{mob_cip}. Since people in group $A$ should only be served by this type of centers, the metric $\mathcal{G}_p$ shows that the population belonging to this group is fully vaccinated. Regarding the participation of the other groups shown with metric $\mathcal{T}_p$, we see that they take advantage of the available capacity of permanent centers and dispense with the use of temporary centers.

\subsection{Management of uncertain demand}
\label{5.7}

The MMV model presented in Section \ref{mathematical_formulation}, defines the vaccination demand in a deterministic way by using the parameter $pop_{lp}$. 
To extend the functionality of the parameter as a part of the model, 
we propose a way to incorporate the uncertainty associated with it.
In the worldwide experience, there are several reasons for this uncertainty, people with a busy schedule not finding time to get vaccinated, lack of access to vaccines, people afraid of vaccines' side effects, people that do not perceive COVID-19 as a threat, people that do not trust in vaccines nor institutions, and a variety of conspiracy theories against vaccines. Also, the level of people's attendance varies according to the neighborhood and the groups defined. Consequently, the (deterministic) uncertainty of the vaccine demand, for each $l \in L_k$ and $p \in P$, could be expressed as follows.

$pop_{lp}= pop^{LB}_{lp}+\eta_{lp}$, where $pop^{LB}_{lp}$ represents a lower bound of the estimation of people in group $p$ to vaccinate in the neighborhood $L$. $\eta_{lp}$ denotes an estimate of the maximum number of additional people that could be vaccinated respect to $pop^{LB}_{lp}$.

To manage the uncertainty, additionally to the parameters $\eta_{lp}$, a new set of variables, $\Gamma_{lp} \in [0,1] $, is defined for each neighborhood, $l$, and group, $p$. The non-robust case is obtained when $\Gamma_{lp}=0 \text{ } \forall l \in L, \forall p \in P$. Each decision variable $\Gamma_{lp}$ greater than zero provides robustness to the model, so it forces to be prepared to attend a bigger demand ($\eta_{lp}\Gamma_{lp}$), so that, a new objective seeks to maximize this robustness. This robustness can be understood like a security inventory of vaccines. Consequently, to mitigate this new objective, the rest of the objective functions are affected. This basic approach to managing the uncertainty allows us to provide a vaccination plan that responds under realistic scenarios. 
 
The changes resulting from the uncertain parameter directly affect the constraints (\ref{r5}) and (\ref{r10}) of the model presented in section \ref{mathematical_formulation}. These constraints are replaced by (\ref{r5.4}) and (\ref{r5.6}), respectively. Constraints (\ref{r5.4}) defines the interval for the vaccine demand. The constraint (\ref{r5.6}) bounds the daily vaccination covering by temporary centers. The objective function (\ref{f_3}) is in charge of managing demand variations.

 


\begin{align}
\sum_{t \in T, i \in I} \phi_{lpti} + \sum_{t \in T, j \in J} \gamma_{lptj} &\geq  pop_{lp} + \eta_{lp}\Gamma_{lp} && \forall l \in \left\lbrace L_{k}: k \in K \right\rbrace, \thinspace p \in P \label{r5.4} \\
v_{jptl} \cdot (pop_{lp}+\eta_{lp}) - \gamma_{lptj}  &\geq 0 && \forall p \in P, \thinspace j \in J, \thinspace l \in \left\lbrace L_{k}: k \in K \right\rbrace \thinspace ,t \in T \label{r5.6}
\end{align}

\begin{align}
\max f_3 = &\sum_{l \in \left\lbrace L_{k}: k \in K \right\rbrace, p \in P} \Gamma_{lp}  \label{f_3}
\end{align}
The model is solved by keeping the weights approach and and can be used as e.g.,
\[ \min Z = \frac{1-\alpha}{2} f_1^\ast + \frac{1-\alpha}{2} f_2^\ast-\alpha f_3^\ast, \]
note that as $f_3$ is seeking to be maximized, it is added to $Z$ with a negative sign. 

\begin{figure}
\centering
\footnotesize
\subfloat []{\includegraphics[width=9cm]{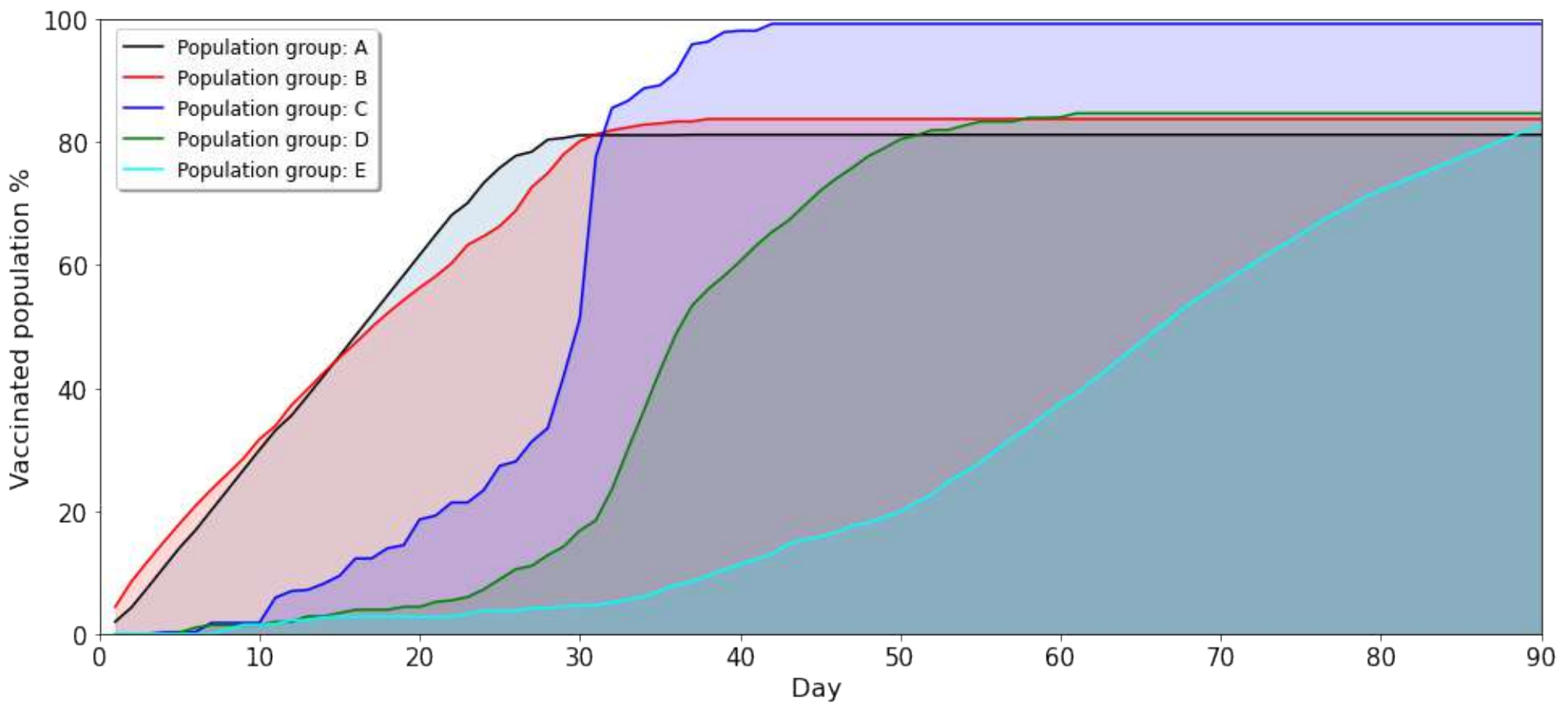}}\hspace{0mm}
\subfloat []{\includegraphics[width=9cm]{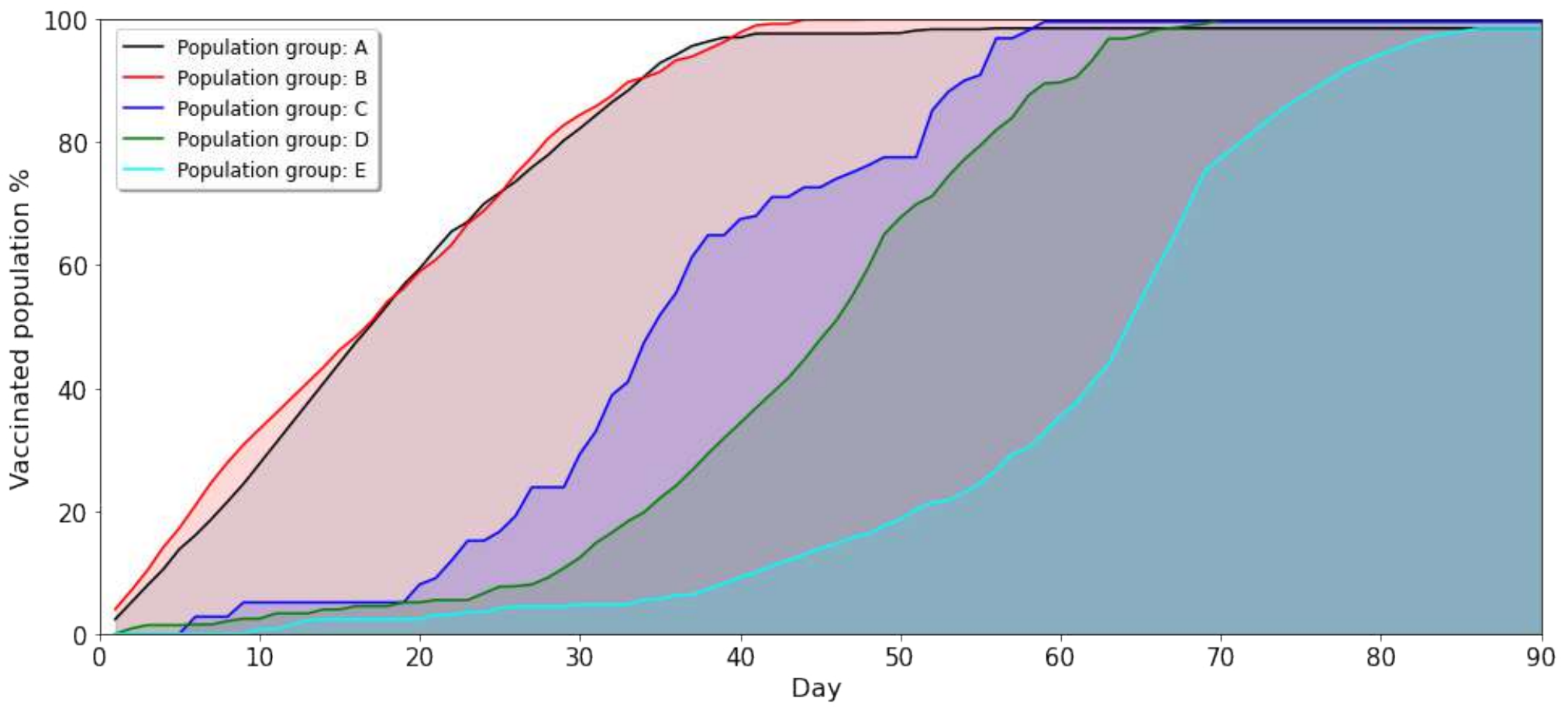}}\vspace{1mm}
\caption{{ \footnotesize  Progress of the vaccination campaign for the different population groups. (a) $\alpha=0.1$ and (b) $\alpha=0.5$
}\label{[plt_incer}}
\label{plt_incer}
\end{figure}

The results of the experimentation of the model using different values of $\alpha$ are presented in Table \ref{gamma}. In this case, $pop^{LB}_{lp}$ is considered to be 80\% of the total number of people to be vaccinated $pop_{lp}$, so $\eta_{lp}$ is responsible for supplying the remaining 20\%. When we look at Table \ref{gamma}, we can see that as we increase the value of the $\alpha$ parameter, the number of people vaccinated increases ($\mathcal{Q}$), even up to the total population($\alpha=0.5$). This increase in the demand for vaccines is accompanied by greater use of temporary care centers ($\mathcal{P}$). Finally, it can be seen that there are small variations in the termination date of the population groups, in comparison with what was presented in previous sections, which can be seen in column $\mathcal{D}_P$ in the Table \ref{gamma} and in the Figures \ref{[plt_incer}-a and b.


\begin{table}
\centering
\begin{tabular}{|c|c|c|c|c|c|c|}
\hline
$\alpha$ & $f_1*$    & $f_2*$    & $f_3*$    & $\mathcal{Q}$ & $\mathcal{P}$&$\mathcal{D}_{A}$-$\mathcal{D}_{B}$-$\mathcal{D}_{C}$-$\mathcal{D}_{D}$-$\mathcal{D}_{E}$\\ \hline
0.1   & 0.34   & 0.17   & 0.44   & 84\% & 38\%&  34-36-41-62-90  \\
0.2   & 0.38   & 0.21   & 0.68   & 89\% & 38\%& 42-88-46-77-90 \\
0.3   & 0.43   & 0.31   & 0.89   & 95\%& 41\% & 38-44-48-80-90\\
0.4   & 0.45 & 0.37 & 0.97 & 97\%& 44\%& 42-44-73-71-90 \\
0.5   & 0.47 & 0.40 & 1.00 & 100\%& 45\%& 56-53-52-71-90 \\ \hline
\end{tabular}
\caption{\footnotesize Results for different $\alpha$ values.$\mathcal{Q}$ represents the percentage represents the percentage of people vaccinated. $\mathcal{P}$ represents the percentage of patients vaccinated in temporary centers. $f_{1}^{*}$, $f_{2}^{*}$ and $f_{3}^{*}$ are the values of the objective functions already normalized. $\mathcal{D}_{p}$ represents the day on which the risk group $p$ $\in$ $\{A,B,C,D,E\}$ has finished being vaccinated. }
\label{gamma}
\end{table}

\subsection{Managerial insights}
\label{5.8}


Important decisions will have to be made in the health sector, specifically regarding vaccination campaigns and lessons learned from the SARS-CoV-2 virus. The 2022 global agenda will be intensely focused on the VSC, seeking to immunize as many people as possible in the shortest time while maintaining social distancing to avoid contagion. Furthermore, considering people with reduced mobility or pre-existing diseases, vaccination plans must consider the use of vaccination centers close to these people. 

Additionally, our model presents reliable support to vaccination campaigns management, scheduling patients' care based on the demand of the various geographical sectors of a given area by means of permanent and temporary centers. Since the application of a vaccine, either against influenza or the COVID-19, is a simple operation, it could be managed in temporary centers like schools and drive-through. Therefore, the use of a significant number of temporary centers in a city is possible, and, consequently, an important part of the population could be vaccinated in these places. 


Through the scenarios defined in the case study, one in which temporary centers have a greater capacity than permanent centers and another in which capacities were equalized, it was possible to obtain the following important findings: 
\begin{itemize}
    \item Firstly, the model's flexibility is highlighted, prioritizing the attention of different population groups under both scenarios.
    \item Secondly, by having vaccination capacity available and playing down the importance of installing temporary centers, the duration of the vaccination campaign can be shortened.  The latter is relevant for decision-makers since it allows them to understand the campaign's length, compared to the costs associated with a less intensive vaccination. 
    \item Finally, the possibility that this model could be used in existing or future pandemics is highlighted.  Furthermore, as a tool for decision-makers, our model would contribute to public sector management, allowing the organization of vaccination campaigns by prioritizing zones and population groups, negotiating the supply of the number of doses, and coordinating with different health actors system.
\end{itemize}

Based on the results previously presented, we highlight the following managerial insights: 
\begin{itemize}
    \item Firstly, an accurate assessment phase for defining the maximum number of temporary centers, $|J|$, must be conducted since these are the main parameter for controlling the congestion. It is clear that the greater the number of temporary centers, the less congestion; so that, before running this model, the optimal number of temporary centers must be defined.  
    \item Note that, this model provides solutions at a daily level; hence, a more detailed plan must be generated once the daily number of persons is defined. A non-informative assumption to estimate the hourly demand can be a constant arrival rate at centers. More sophisticated assumptions based on data or expert knowledge can be considered to model this arrival rate. 
    \item \ To deploy the solution from our model, it should be implemented using a computational platform that generates both the daily schedule (tactical plan) and the hourly arrival for each day and mobile team (operative plan). Additionally, this system should consider a channel to interact with users (i.e., App or web page). Current vaccination campaigns in Chile use a computer platform called the National Immunization Registry, which allows health care workers to record information about vaccinated people to control the vaccination process. The system is used throughout the whole country. It serves to monitor the progress of a campaign virtually and control the stock of vaccines. Our model should be incorporated into a system like this.  
    \item  The number of temporary centers in each neighborhood in the baseline model is assumed to be equal to one; however, we recommend adapting (increase) this parameter when the scenario includes highly heterogeneous neighborhoods. This setting is relevant when you have crowded communities covering extensive geographic zones.
    \item  Depending on the scenario, the planner must select an appropriate combination of weights and the number of objectives. For example, if there is no reliable system to manage the vaccine campaign, we recommend setting high importance to the weight associated with uncertainty.
\end{itemize}


\section{Conclusions}
\label{section6}

This paper proposed a bi-objective optimization model to plan a vaccination campaign during the COVID-19 pandemic. Planning a campaign presents critical coordination challenges between different actors in the health system. New requirements were added, such as avoiding congestion and recognizing that some people are more vulnerable to the virus than others. In response to this, our model considers the use of permanent and temporary centers, allowing the care of elderly patients with mobility problems while managing crowds avoiding congestion. 


The baseline model has two objectives. Firstly, by optimizing an ad-hoc priority function, the model focuses on executing a vaccination plan as soon as possible, prioritizing high-risk groups. Secondly, an economic objective was included that is associated with the deployment of temporary centers. From the trade-off between these objectives, it was possible to obtain vaccination plans that meet the decision-makers requirements. In order to manage uncertainty in the daily arrives of people, we designed a reformulated model with an additional objective. 

To test the model, a case based on single-dose vaccination campaign was analyzed, and the study area was the commune of San Bernardo, Chile. In total, 115,800 doses were administered to different population groups, where depending on the vaccination centers' capacities, the campaign's duration could be reduced by up to 27\%. Because vaccination campaigns related to influenza or COVID-19 are being executed worldwide, the proposed model could be used in different locations, providing the option of adapting it to existing or future pandemics. Additionally, using our model, decision-makers can have a more holistic view of planning vaccination campaigns, including important insights that can be available to various actors in the VSC. Finally, integrating information systems and mathematical programming models, such as the one presented in this work, can provide real-time information to generate more robust plans.

Based on our model, a robust variant should be formulated for future research, including uncertain supply and vaccines' reliability. From an algorithmic perspective, efficient heuristics can be designed when challenging size problems are faced. We foresee that the challenge of better integrating production and distribution will remain critical, especially in developing countries. Therefore, new stages could be added to the model and further explore the results.

\bibliographystyle{unsrt}  
\bibliography{references}

\end{document}